\documentclass[12pt]{article}
\usepackage[cp1251]{inputenc}
\usepackage[russian]{babel}
\usepackage{amssymb}
\usepackage{amsmath}

\newcommand{\fu}{\buildrel{_\wedge}\over}

\newcommand{\rightfu}{{^\wedge}}
\newcommand{\rightuf}{{^\vee}}
\newcommand{\widefu}{\widehat}

\begin{document}

\begin{center}
\textbf{\large Quantitative estimates in Beurling--Helson type
theorems}\footnote{Published in \emph{Sbornik: Mathematics},
\textbf{201}:12 (2010), 1811-1836. The text below may slightly
vary from the finally printed version.}
\end{center}

\begin{center}
 Vladimir Lebedev
\end{center}

\begin{quotation}
{\small \textsc{Abstract.} We consider the spaces $A_p(\mathbb T)$
of functions $f$ on the circle $\mathbb T$ such that the sequence
of Fourier coefficients $\fu{\!f}=\{\fu{\!f}(k), ~k \in \mathbb
Z\}$ belongs to $l^p, ~1\leq p<2$. The norm on $A_p(\mathbb T)$ is
defined by $\|f\|_{A_p}=\|\fu{\!f}\nolinebreak\|_{l^p}$. We study
the rate of growth of the norms $\|e^{i\lambda\varphi}\|_{A_p}$ as
$|\lambda|\rightarrow \infty, ~\lambda\in\mathbb R,$ for $C^1$
-smooth real functions $\varphi$ on $\mathbb T$. The results have
natural applications to the problem on changes of variable in the
spaces $A_p(\mathbb T)$.

  References: 17 items.

  }

AMS 2010 Mathematics Subject Classification 42A16
\end{quotation}

\quad

\begin{center}
\textbf{Introduction}
\end{center}

    We consider Fourier series
$$
f(t)\sim \sum_{k\in \mathbb{Z}}\fu{\!f}(k)e^{ikt}
$$
of (integrable) functions $f$ on the circle
$\mathbb{T=R}/2\pi\mathbb Z$ ($\mathbb R$ is the real line,
$\mathbb Z$ is the set of integers).

   Let $A_1(\mathbb T)$ be the space of continuous functions
$f$ on $\mathbb T$ such that the sequence of Fourier coefficients
$\fu{\!f}=\{\fu{\!f}(k), ~k\in{\mathbb Z}\}$ belongs to $l^1$. Let
$A_p(\mathbb T),\\ ~1<p\leq 2,$ be the space of integrable
functions $f$ on $\mathbb T$ such that $\fu{\!f}$ belongs to
$l^p$. Endowed with the natural norms
$$
\|f\|_{A_p(\mathbb T)}=\|\fu{\!f}\|_{l_p}=\bigg(\sum_{k\in \mathbb{Z}}
|\fu{\!f}(k)|^p\bigg)^{1/p},
$$
the spaces $A_p(\mathbb T), ~1\leq p\leq 2,$ are Banach spaces.
The space $A(\mathbb T)=A_1(\mathbb T)$ is a Banach algebra (with
the usual multiplication of functions).

   Suppose that we have a map of the circle into itself,
i.e. a function $\varphi : \mathbb R\rightarrow\mathbb R$ such
that $\varphi(t+2\pi)=\varphi(t) ~(\textrm{mod} ~2\pi)$. According
to the Beurling--Helson theorem [1] (see also [2], [3]), if
$\|e^{in\varphi}\|_{A(\mathbb T)}=O(1), ~n \in \mathbb Z,$ then
the function $\varphi$ is $\mathrm{mod}\, 2\pi$ -linear (affine)
with integer tangent coefficient: $\varphi(t)=mt+t_0 (\mathrm{mod
\, 2\pi}), ~m\in \mathbb Z$. This theorem yields the solution of
the Levy problem on the description of endomorphisms of the
algebra $A(\mathbb T)$; all these endomorphisms prove to be
trivial, i.e., have the form $f(t)\rightarrow f(mt+t_0)$. In other
words only trivial changes of variable are admissible in
$A(\mathbb T)$. We note also another version of the
Beurling--Helson theorem: if $U$ is a bounded
translation-invariant operator from $l^1$ to itself such that
$\|U^n\|_{l^1\rightarrow l^1}=O(1), ~n \in \mathbb Z$, then $U=\xi
S$, where $\xi$ is a constant, $|\xi|=1$, and $S$ is a
translation.

  At the same time, although the Beurling--Helson theorem
states unbounded growth of the norms $\|e^{in\varphi}\|_A$ for
nonlinear maps $\varphi: \mathbb T \rightarrow \mathbb T$, the
rate of growth of these norms as $|n|\rightarrow\infty$ is in
general not clear. The same concerns the behavior of the norms
$\|e^{in\varphi}\|_{A_p}, ~p>1$.

  Let $C^\nu(\mathbb T)$ be the class of (complex-valued)
functions on $\mathbb T$ with continuous derivative of order
$\nu$. We have $C^1(\mathbb T)\subseteq A(\mathbb T)\subseteq
A_p(\mathbb T)$.

   In this paper we mainly consider real
$C^1$ -smooth functions $\varphi$ on $\mathbb T$ and study the
growth of the norms $\|e^{i\lambda\varphi}\|_{A_p}$ as
$|\lambda|\rightarrow\infty, ~\lambda\in\mathbb R$. The only
interesting case is certainly that of non-constant functions
$\varphi$. The corresponding results on the behavior of the
exponential functions $e^{in\varphi}$ for nonlinear maps $\varphi:
\mathbb T\rightarrow\mathbb T$ and integer frequencies $n$ will be
obtained as simple corollaries.

  It is not difficult to show that if $\varphi\in C^1(\mathbb T)$
is a real function (and moreover if $\varphi$ is an absolutely
continuous real function  with the derivative in $L^2(\mathbb T)$)
then for $1\leq p<2$ we have
$$
\|e^{i\lambda \varphi}\|_{A_p(\mathbb T)}=
O(|\lambda|^{\frac{1}{p}-\frac{1}{2}}),
\qquad |\lambda|\rightarrow\infty, \quad\lambda\in\mathbb R,
\eqno(1)
$$
(see [2, Ch. VI, \S ~3] in the case when $p=1$; for $1<p<2$ the
estimate follows immediately by interpolation between $l^1$ and
$l^2$).

   On the other hand, lower estimates for the norms
of $e^{i\lambda\varphi}$ in $A_p$ for functions $\varphi$ of class
$C^2$ have long been known. Suppose that $\varphi\in C^2(\mathbb
T)$ is a real non-constant function and $1\leq p < 2$. Then
$$
\|e^{i\lambda\varphi}\|_{A_p(\mathbb T)}\geq c\, |\lambda|^{
\frac{1}{p}-\frac{1}{2}}, \qquad \lambda \in \mathbb R, \eqno(2)
$$
where $c=c(p,\varphi)>0$ is independent of $\lambda$. For $p=1$
this estimate is contained implicitly in Leibenson's paper [4] and
was obtained in explicit form by Kahane [5] who used Leibenson's
approach. In the general case estimate (2) was obtained (by the
same method) by Alpar [6]. A short and simple proof for $p=1$ can
be found in [2, Ch. VI, \S~3] and in the general case in [7].

  Thus, if $\varphi \in C^2(\mathbb T)$ is a real function,
$\varphi\neq\mathrm{const}$, then
$$
\|e^{i\lambda\varphi}\|_{A_p(\mathbb T)}\simeq
|\lambda|^{\frac{1}{p}-\frac{1}{2}}
$$
for all $p, ~1\leq p<2$. In particular
$\|e^{i\lambda\varphi}\|_A\simeq \sqrt{|\lambda|}$. (The sign
$\simeq$ means that for all sufficiently large $|\lambda|$ the
ratio of the corresponding quantities is contained between two
positive constants).

    We note that the proof of the Leibenson--Kahane--Alpar estimate
(2) is based on the second part of the van der Corput lemma (see
[8, Ch. V, Lemma 4.3]) and essentially uses the condition that the
curvature of some ark of the graph of $\varphi$ is bounded away
from zero, i.e., the condition $|\varphi ''(t)| \geq \rho
>0, ~t \in I,$ where $I$ is a certain interval.

  In general, the norms $\|e^{i\lambda\varphi}\|_A$ can grow
rather slowly. Kahane showed (see [2, Ch. VI, \S ~2]) that if
$\varphi\neq\mathrm{const}$ is a real continuous piecewise linear
function on $\mathbb T$ (which means that $[0, 2\pi]$ is a finite
union of some intervals such that $\varphi$ is linear on each of
them) then $\|e^{i\lambda\varphi}\|_A\simeq \log |\lambda|$. It is
not clear if there are non-trivial functions $\varphi$ that yield
the growth of the norms of the exponential functions in $A(\mathbb
T)$ slower than logarithmic. Kahane conjectured [2], [3] that the
Beurling--Helson theorem can be strengthened considerably. In
particular, he posed the following question: is it true that if
$\|e^{i\lambda\varphi}\|_A=o(\log |\lambda|)$, then
$\varphi=\mathrm{const}$. We do not know the answer even under the
assumption that $\varphi\in C^1(\mathbb T)$.

  We also note that for any piecewise linear
real function $\varphi$ on $\mathbb T$ we have
$\|e^{i\lambda\varphi}\|_{A_p}=O(1)$ for all $p>1$ (see  [7]).

   We recall the results known for
the $C^1$ -smooth case (besides estimate~(1)).

   In [7] (a joint work of the author and Olevski\v{\i})
we constructed a real function $\varphi\in C^1(\mathbb
T),~\varphi\neq\mathrm{const},$ such that
$\|e^{i\lambda\varphi}\|_{A_p}=O(1)$ for all $p>1$. In addition
this function is nowhere linear, i.e., it is not linear on any
interval (and thus, in a sense, it is essentially different from
piecewise linear functions).

   In [9] the author of the present paper showed
that for functions $\varphi \in C^1$ the norms
$\|e^{i\lambda\varphi}\|_A$ can grow rather slowly, namely: if
$\gamma (\lambda)\geq 0$ and $\gamma (\lambda) \rightarrow
+\infty$ as $\lambda\rightarrow+\infty$, then there exists a
nowhere linear real function  $\varphi \in C^1 (\mathbb T)$ such
that
$$
\|e^{i\lambda\varphi}\|_A=O(\gamma (|\lambda|)\log
|\lambda|).
$$

  Thus, the case of $C^1$ -smooth phase $\varphi$ is essentially
different from the $C^2$ -smooth case.

  As far as we know, the only lower estimate for
the norms $\|e^{i\lambda\varphi}\|_A$ obtained previously in the
case when $\varphi\in C^1$ but the twice differentiability of
$\varphi$ is not assumed is due to Leblanc [10]: if a real
function $\varphi\in C^1(\mathbb T)$ is non-constant and its
derivative $\varphi'$ satisfies the Lipschitz condition of order
 $\alpha, ~0<\alpha\leq 1,$ then
$$
\|e^{i\lambda\varphi}\|_A\geq c
\frac{|\lambda|^{\frac{\alpha}{1+\alpha}}}{(\log |\lambda|)^2}, \qquad
|\lambda|\geq 2. \eqno(3)
$$

   Now we shall briefly describe the results of the present work.
In \S ~1 we prove Theorem 1 in which we give lower estimates for
$\|e^{i\lambda\varphi}\|_{A_p}$ for $C^1$ -smooth real functions
$\varphi$ on $\mathbb T$. In \S ~2 we obtain Theorem 2, which
shows that the estimates of Theorem 1 are close to being sharp and
in certain cases are sharp. Here, in the introduction, we omit the
statements of the theorems and only note their corollaries.

  From Theorem 1 it follows that if a real function
$\varphi\in C^1(\mathbb T)$ is non-constant and $\varphi'$
satisfies the Lipschitz condition of order $\alpha, ~0<\alpha \leq
1,$ then
$$
\|e^{i\lambda\varphi}\|_{A_p}\geq c_p
|\lambda|^{\frac{1}{p}-\frac{1}{1+\alpha}}, \quad \lambda\in\mathbb R,
\eqno(4)
$$
for all $p, ~1\leq p<1+\alpha$ (see Corollary 1). In particular,
putting $p=1$ here, we obtain the result, which is stronger than
Leblanc's estimate (3), namely,
$$
\|e^{i\lambda\varphi}\|_A \geq c
|\lambda|^{\frac{\alpha}{1+\alpha}},\quad \lambda\in\mathbb R.
$$
For $\varphi\in C^2$ we have $\alpha=1$ and from (4) we get the
Leibenson--Kahane--Alpar estimate (2).

  We note that the Leibenson--Kahane--Alpar estimate and
the Leblanc estimate have local character; roughly speaking, they
remain valid if we assume that $\varphi$ is nonlinear on some
interval and has an appropriate smoothness on this interval. Our
lower estimates are of local character as well (see Theorem 1$'$).

   A particular case of Theorem 2 is the following statement:
for each $\alpha, ~0<\alpha<1,$ there exists a nowhere linear real
function $\varphi\in C^1(\mathbb T)$ such that its derivative
$\varphi'$ satisfies the Lipschitz condition of order $\alpha$ and
we have
$$
\|e^{i\lambda \varphi}\|_A=O\big(|\lambda|^{\frac{\alpha}{1+\alpha}}
(\log |\lambda|)^{\frac{1-\alpha}{1+\alpha}}\big), \qquad
|\lambda|\rightarrow\infty
$$
(see Corollary 2). For the same function $\varphi$ we have
$$
\|e^{i\lambda \varphi}\|_{A_p}\simeq
|\lambda|^{\frac{1}{p}-\frac{1}{1+\alpha}}
\quad \textrm{if} \quad 1<p<1+\alpha
$$
 and
$$
\|e^{i\lambda \varphi}\|_{A_p}\simeq 1
\quad \textrm{if} \quad 1+\alpha<p<2.
$$
In addition,
$$
\|e^{i\lambda \varphi}\|_{A_p}=O((\log |\lambda|)^{1/p})
\quad \textrm{for} \quad p=1+\alpha.
$$

  From estimate (4) and the obvious estimate
$$
\|e^{i\lambda \varphi}\|_{A_p(\mathbb T)}\geq \|e^{i\lambda
\varphi}\|_{A_2(\mathbb T)}=\|e^{i\lambda \varphi}\|_{L^2(\mathbb T)}=1,
\qquad 1\leq p\leq 2, \eqno(5)
$$
it follows that the function $\varphi$ that we have constructed
yields the slowest possible growth of the norms $\|e^{i\lambda
\varphi}\|_{A_p}$ for $1<p<2, ~p\neq 1+\alpha$.

  As we mentioned above, there exists a non-trivial $C^1$ -smooth
function $\varphi$ with extremely slow (arbitrarily close to
logarithmic) growth of the norms $\|e^{i\lambda\varphi}\|_A$. This
result, which the author proved earlier, follows immediately from
Theorem 2 (see Corollary 3).

   In \S ~3 we consider $C^1$ -smooth maps of the circle
into itself and give the corresponding versions of the results
obtained in \S\S 1 and 2. These versions have natural applications
in the study of change of variable operators (super-\\position
operators) $f\rightarrow f\circ \varphi$ on the spaces $A_p$. In
particular, we show how smooth a nonlinear map $\varphi : \mathbb
T\rightarrow \mathbb T$ can be provided that $f\circ
\varphi\in\bigcap_{p>1} A_p$ whenever $f\in A$.

\smallskip

  By $\omega (I, g, \delta )$ we denote the modulus of continuity
of a function $g$ on an interval $I\subseteq\mathbb R$:
$$
\omega (I,g,\delta )=\sup_{|t_1-t_2|\leq\delta \atop {t_1, t_2 \in
I}}|g(t_1)-g(t_2)|, \qquad \delta\geq 0.
$$
If $I=\mathbb R$, then we just write $\omega (g,\delta )$. The
class $\mathrm{Lip}_\omega (I)$ consists of all functions $g$ on
$I$ with $\omega (I,g,\delta )=O(\omega (\delta )), ~\delta \to
+0$, where $\omega (\delta)$ is a given continuous non-decreasing
function on $[0,+\infty),~ \omega (0)=0$. The class $C^{1,
\omega}(I)$ consists of functions $g$ on $I$ with the derivative
$g'\in \mathrm{Lip}_\omega (I)$. The class $C^{1, \omega}(\mathbb
T)$ consists of $2\pi$ -periodic functions that belong to $C^{1,
\omega}(\mathbb R)$. For $0<\alpha \leq 1$ we write $C^{1,
\alpha}$ instead of $C^{1, \delta ^\alpha}$. The class $C^1(I)$
consists of continuously differentiable functions on $I$. For a
set $E\subseteq\mathbb T$ we denote its characteristic function by
$1_E$: $1_E(t)=1$ for $t\in E$, $1_E(t)=0$ for $t\in \mathbb
T\setminus E$ and (if $E$ is measurable) by $|E|$ its (Lebesgue)
measure. If $f$ and $E$ are a function on $\mathbb T$ and a set in
$\mathbb T$, respectively, we write $\mathrm{supp}\,f\subseteq E$
if $f(t)=0$ for almost all $t\in \mathbb T\setminus E$. We use the
same notation for sets in $[0, 2\pi]$ or in $\mathbb R$ and for
functions defined on $[0, 2\pi]$ or on $\mathbb R$. In the
standard way we identify integrable functions on $[0, 2\pi]$ with
integrable functions on $\mathbb T$. We use $c, c_1, c_p, c_{p,
1}, c(p)$, etc. to denote various positive constants which may
depend only on $p$ and $\varphi$.

\quad

\begin{center}
\textbf{\S~1. Lower estimates}
\end{center}

   \textbf{Theorem 1.} \emph {Let $1\leq p<2$. Let $\varphi$ be a
real function on $\mathbb T$. Suppose that $\varphi$ is
non-constant and $\varphi \in C^{1, \omega}(\mathbb T)$. Then
$$
\|e^{i\lambda\varphi}\|_{A_p(\mathbb T)}\geq c\, |\lambda|^{1/p} \,\chi
^{-1} (1/|\lambda|), \qquad \lambda \in \mathbb R, \quad |\lambda|\geq 1,
$$
where $\chi ^{-1}$ is the function inverse to $\chi (\delta
)=\delta \omega (\delta )$ and $c=c(p, \varphi)>0$ is independent
of $\lambda$.}

\quad

  Certainly in view of the obvious estimate (5),
the estimate in Theorem 1 makes sense only in the case when its
right-hand side is growing unboundedly together with $\lambda$.
(This is always the case when $p=1$).

    Theorem 1 immediately implies the following corollary.

\quad

   \textbf{Corollary 1.} \emph {Let $0<\alpha\leq 1$.
If $\varphi$ is a real non-constant function on $\mathbb T$ and
$\varphi \in C^{1, \alpha}(\mathbb T)$ then for all $p, ~1\leq
p<1+\alpha,$ we have
$$
\|e^{i\lambda\varphi}\|_{A_p(\mathbb T)}\geq c_p\,
|\lambda|^{\frac{1}{p}-\frac{1}{1+\alpha}}, \qquad \lambda \in \mathbb R.
$$
In particular, $\|e^{i\lambda\varphi}\|_{A(\mathbb T)}\geq c\,
|\lambda|^{\frac{\alpha }{1+\alpha}}$.}

\quad

   As we mentioned in Introduction, a particular case of Corollary 1
is the Leibenson--Kahane--Alpar estimate. For $p=1$ the estimate
of Corollary 1 improves the result of Leblanc.

  Theorem 1$'$ below is a local version of Theorem 1.

  Let $I\subset\mathbb R$ be an interval of length less then $2\pi$.
We say that a function $f$ defined on $I$ belongs to $A_p(\mathbb
T, I)$ if there is a function $F$ in $A_p(\mathbb T)$ such that
its restriction $F_{|I}$ to $I$ coincides with $f$. We put
$$
\|f\|_{A_p(\mathbb T, I)}=
\inf_{F\in A_p(\mathbb T), ~F_{|I}=f} \|F\|_{A_p(\mathbb T)}.
$$
It is clear that $A_p(\mathbb T, I)$ is a Banach space, $1\leq
p\leq 2$.

\quad

     \textbf{Theorem 1$'$.} \emph {Let $1\leq p<2$. Let
$\varphi$ be a real function on an interval $I\subset\mathbb R,
~|I|<2\pi$. Suppose that $\varphi$ is nonlinear on $I$ and
$\varphi \in C^{1, \omega}(I)$. Then}
$$
\|e^{i\lambda\varphi}\|_{A_p(\mathbb T, I)}
\geq c\, |\lambda|^{1/p} \,\chi ^{-1}(1/|\lambda|),
\qquad \lambda \in \mathbb R, \quad |\lambda|\geq 1.
$$

\quad

   The local version of Corollary 1 is obvious.

\quad

    \emph{Proof of Theorem} 1. Let
$$
m=\min_{t \in [0, 2\pi]} \varphi'(t), \quad M=\max_{t \in [0,
2\pi]} \varphi'(t).
$$
Since $\varphi$ is non-constant whereas every continuous $2\pi$
-periodic function linear on $[0, 2\pi]$ is constant, we have
$m<M$.

   Fix $c>0$ so that
$$
\omega (\varphi', \delta)\leq c \omega (\delta), \quad \delta\geq 0.
$$

  Let $\lambda \in \mathbb R$. Without loss of generality we can
assume that $\lambda>0$ (complex conjugation does not affect the
norm of a function in $A_p$). Define $\delta _\lambda>0$ by
$$
\chi (2\delta_\lambda)=\frac{1}{2c\lambda}. \eqno(6)
$$

   For $0<\varepsilon \leq \pi $ let
$\Delta _\varepsilon$ be the ``triangle'' function supported on
the interval $(-\varepsilon, \varepsilon)$, i.e. the function on
$\mathbb T$, defined as follows:
$$
\Delta _\varepsilon (t)=
\max \bigg (1-\frac{|t|}{\varepsilon },~0 \bigg),
\qquad t \in [-\pi ,~\pi].
$$

  For an arbitrary interval $J\subseteq [0, 2\pi]$ let $\Delta _J$
be the triangle function supported on $J$, namely:
$\Delta_J(t)=\Delta _{|J|/2}(t-c_J)$, where $c_J$ is the center of
$J$ (and $|J|$ is its length).

\quad

\textbf{Lemma 1.} \emph{Suppose that $\lambda>0$ is sufficiently
large. Then for every $k\in \mathbb Z$ satisfying
$m\lambda<k<M\lambda$ there exists an interval $I_{\lambda,
k}\subseteq [0, 2\pi]$ such that $|I_{\lambda,
k}|=2\delta_\lambda$ and
$$
|(\Delta _{I_{\lambda, k}} e^{i\lambda \varphi})^\rightfu (k)|\geq
\frac{\delta_\lambda}{4\pi}.
$$
}

\quad

\emph{Proof.} We assume $\lambda$ to be so large, that
$$
2\delta_\lambda<2\pi \eqno(7)
$$
and
$$
\lambda >\frac {2}{M-m}. \eqno(8)
$$

   Take an arbitrary $k \in \mathbb Z$ with
$m\lambda < k < M\lambda$ (see (8)). We can find a point
$t_{\lambda, k}$ in $(0, 2\pi)$ such that $\varphi'(t_{\lambda,
k})= k/\lambda$. There exists an interval $I_{\lambda, k}\subseteq
[0, 2\pi]$ that contains the point $t_{\lambda, k}$ and has length
$2\delta_\lambda$ (see (7)). Consider the following linear
function:
$$
\varphi_{\lambda, k}(t)= \varphi(t_{\lambda, k})+\frac{k}{\lambda}
(t-t_{\lambda, k}), \qquad t \in [0, 2\pi].
$$

   If $t \in I_{\lambda, k}$ then for some point $\theta$
that lies between $t$ and $t_{\lambda, k}$ we have
$$
\varphi(t)-\varphi(t_{\lambda, k})=\varphi'(\theta ) (t-t_{\lambda, k}).
$$
So,
$$
|\varphi(t)-\varphi_{\lambda, k}(t)|=|(\varphi(t)-\varphi(t_{\lambda,
k}))-\frac{k}{\lambda} (t-t_{\lambda, k})|
$$
$$
=|\varphi'(\theta ) (t-t_{\lambda, k})-\varphi'(t_{\lambda, k})
(t-t_{\lambda, k})|
$$
$$
=|t-t_{\lambda, k}| |\varphi'(\theta)-\varphi'(t_{\lambda, k})| \leq
2\delta _\lambda \omega (\varphi', 2\delta _\lambda)\leq 2\delta _\lambda
c\omega (2\delta _\lambda)=c\chi (2\delta _\lambda).
$$
Hence, taking into account (6), we have
$$
|e^{i\lambda\varphi(t)}-e^{i\lambda\varphi_{\lambda, k}(t)}| \leq
|\lambda\varphi(t)-\lambda\varphi_{\lambda, k}(t)|\leq \lambda c \chi
(2\delta _\lambda)=\frac{1}{2}, \qquad t \in I_{\lambda, k}. \eqno(9)
$$

     Using (9), we obtain
$$
|(\Delta _{I_{\lambda, k}} e^{i\lambda\varphi})^\rightfu(k)-(\Delta
_{I_{\lambda, k}} e^{i\lambda\varphi_{\lambda, k}})^\rightfu(k)| \leq
\frac{1}{2\pi}\int_0^{2\pi} \Delta _{I_{\lambda, k}}(t)
|e^{i\lambda\varphi(t)}-e^{i\lambda\varphi_{\lambda, k}(t)}|~dt
$$
$$
\leq \frac{1}{2}\cdot\frac{1}{2\pi} \int_0^{2\pi} \Delta _{I_{\lambda,
k}}(t)~dt =\frac {1}{2}\widefu{\Delta_{I_{\lambda, k}}}(0).
$$
At the same time,
$$
|(\Delta _{I_{\lambda, k}} e^{i\lambda\varphi_{\lambda, k}})^\rightfu
(k)|=\bigg| \frac{1}{2\pi}\int_0^{2\pi} \Delta _{I_{\lambda, k}}(t)
e^{i(\lambda\varphi_{\lambda, k}(t)-kt)} ~dt \bigg|=
$$
$$
=\frac{1}{2\pi}\int_0^{2\pi} \Delta _{I_{\lambda, k}}(t)~dt=\widefu
{\Delta _{I_{\lambda, k}}}(0),
$$
and thus
$$
|(\Delta _{I_{\lambda, k}} e^{i\lambda \varphi})^\rightfu (k)|
\geq \frac{1}{2}\widefu {\Delta _{I_{\lambda, k}}}(0)=\frac
{\delta_\lambda}{4\pi}.
$$
The lemma is proved.

\quad

  For each $\lambda\in\mathbb R$ define a function $g_\lambda$
by its Fourier expansion:
$$
g_\lambda(t)\sim \sum_{k \in \mathbb Z}
|\widefu{e^{i\lambda\varphi}}(k)| \, e^{ikt}.
$$
Obviously, $g_\lambda\in L^2(\mathbb T)$.

  It is well known that $\Delta_\varepsilon\in A(\mathbb T)$ and
the Fourier coefficients of $\Delta_\varepsilon$ are
non-negative.\footnote{Direct calculation yields $\widefu
{\Delta_\varepsilon}(k)= \frac {2}{\pi}\frac{\sin^2 (\varepsilon
k/2)}{\varepsilon k^2}, ~k\neq 0, \quad \widefu
{\Delta_\varepsilon }(0)=\frac{\varepsilon}{2\pi}$.} Since for an
arbitrary interval $J$ the function $\Delta_{J}$ is obtained from
$\Delta_{|J|/2}$ by translation, we have
$|\widefu{\Delta_{J}}(k)|=\widefu{\Delta_{|J|/2}}(k), ~k \in
\mathbb Z$. Thus, from Lemma 1 it follows that if $\lambda$ is
sufficiently large, then
$$
\frac{\delta _\lambda}{4\pi}\leq |(\Delta _{I_{\lambda, k}} e^{i\lambda
\varphi})^\rightfu (k)|
$$
$$
=\bigg| \sum_{\nu \in \mathbb Z} \widefu{\Delta _{I_{\lambda, k}}}(\nu)
\widefu{e^{i\lambda\varphi }}(k-\nu ) \bigg| \leq \sum_{\nu \in \mathbb
Z} |\widefu{\Delta_{I_{\lambda, k}}}(\nu
)||\widefu{e^{i\lambda\varphi}}(k-\nu )|
$$
$$
=\sum_{\nu \in \mathbb Z}\widefu{\Delta_{\delta_\lambda}}(\nu) \widefu
{g_\lambda}(k-\nu ) = (\Delta _{\delta_\lambda}
\cdot g_\lambda)^\rightfu (k),
\qquad m\lambda< k< M\lambda. \eqno(10)
$$

  Since $\|\Delta _\varepsilon \|_A=\sum_k \widefu
{\Delta_\varepsilon }(k)=\Delta _\varepsilon (0)=1$, it follows
that for every function $f\in A_p(\mathbb T)$ we have
$\|\Delta_\varepsilon \cdot f\|_{A_p}\leq \|f\|_{A_p}, ~1\leq p
\leq 2, ~0<\varepsilon \leq \pi$. So, using (10), we see that for
all sufficiently large $\lambda$
$$
\bigg (\frac{1}{2}(M-m) \lambda \bigg )^{1/p} \,\frac{\delta
_\lambda}{4\pi} \leq \bigg(\sum_{m\lambda<k<M\lambda} \bigg(\frac{\delta
_\lambda}{4\pi} \bigg)^p \bigg)^{1/p}
$$
$$
\leq \bigg(\sum_{m\lambda<k<M\lambda} |(\Delta _{\delta_\lambda} \cdot
g_\lambda)^\rightfu (k)|^p \bigg)^{1/p}
$$
$$
\leq\|\Delta_{\delta _\lambda} \cdot g_\lambda \|_{A_p} \leq
\|g_\lambda\|_{A_p}=\|e^{i\lambda\varphi}\|_{A_p}. \eqno(11)
$$

   Note now that condition (6) yields
$$
1=2c \lambda 2\delta_\lambda \omega(2\delta_\lambda)\leq
4c\lambda\delta_\lambda \omega(4\delta_\lambda)\leq \lambda (c+1)
4\delta_\lambda\omega
((c+1)4\delta_\lambda)=\lambda\chi((c+1)4\delta_\lambda),
$$
and therefore,
$$
\delta_\lambda\geq \frac{1}{4(c+1)}\chi ^{-1}(1/\lambda).
$$
Substituting this estimate in (11) we obtain the statement of the
theorem.

\quad

  \emph{Proof of Theorem} 1$'$. We can obtain it by
an obvious modification of the proof of Theorem 1. Namely,
assuming that $I$ is a closed subinterval of $[0, 2\pi]$ (this
does not restrict generality), we put $m=\min_{t \in I}
\varphi'(t)$ and $M=\max_{t \in I} \varphi'(t)$. Since $\varphi$
is nonlinear on $I$, we have $m<M$. Fix $p, ~1\leq p<2$. For every
$\lambda$ fix a function $F_\lambda\in A_p(\mathbb T)$ which is a
$2\pi$ -periodic extension of $e^{i\lambda\varphi}$ from $I$ to
$\mathbb R$ such that
$$
\|F_\lambda\|_{A_p(\mathbb T)}\leq
2\|e^{i\lambda\varphi}\|_{A_p(\mathbb T, I)}.
$$
In the proof of Theorem 1 one should replace $\omega (\varphi',
2\delta_\lambda)$ by $\omega (I, \varphi', 2\delta_\lambda)$ and,
assuming that $\lambda$ is sufficiently large, one should replace
relation (7) by $2\delta_\lambda<|I|$. For $m\lambda<k<M\lambda$
we choose a point $t_{\lambda, k}$ with $\varphi'(t_{\lambda,
k})=k/\lambda$ so that it lies in the interior of $I$. We choose
an interval $I_{\lambda, k}$ of length $2\delta_\lambda$, which
contains $t_{\lambda, k}$, so that it lies in $I$. Instead of
$e^{i\lambda\varphi}$ one should consider $F_\lambda$.

\quad

\emph{Remark} 1. As we mentioned in Introduction, we have
$$
\|e^{i\lambda\varphi}\|_{A_p(\mathbb T)}=O\big(|\lambda|^{\frac{1}{p}-
\frac{1}{2}}\big), \qquad |\lambda|\rightarrow\infty,
$$
for any real function $\varphi\in C^1(\mathbb T)$. So from
Corollary 1 we see that for any real function $\varphi\in C^{1,
1}(\mathbb T), ~\varphi\neq\mathrm{const},$ we have
$$
\|e^{i\lambda\varphi}\|_{A_p(\mathbb
T)}\simeq|\lambda|^{\frac{1}{p}- \frac{1}{2}}.
$$
Thus, the rate of growth of the norms
$\|e^{i\lambda\varphi}\|_{A_p}$ in the case of a $C^{1, 1}$
-smooth phase $\varphi$ is the same as in the $C^2$ -smooth case.

\quad

  \emph{Remark} 2. A simple modification of the
proof of Theorems 1, 1$'$ allows to obtain their multidimensional
versions. Our results on behavior of the norms
$\|e^{i\lambda\varphi}\|_{A_p(\mathbb T^m)}$ for smooth functions
$\varphi$ on the torus $\mathbb T^m, ~m\geq 2,$ will be presented
elsewhere.

\quad

\begin{center}
\textbf{\S~2. Slow growth of $\|e^{i\lambda\varphi}\|_{A_p(\mathbb
T)}$}
\end{center}

  In this section, for each class $C^{1, \omega}$ (under a certain simple
condition imposed on $\omega$), we shall construct a non-trivial
function $\varphi\in C^{1, \omega}(\mathbb T)$ such that the norms
$\|e^{i\lambda\varphi}\|_{A_p(\mathbb T)}$ have slow growth.

  The case when $\omega(\delta)=\delta$ is described in Remark 1.

  Recall that we say that a function is
nowhere linear if it is not linear on any interval.

  As above, $\chi^{-1}$ is the function inverse to
$\chi(\delta)=\delta\omega(\delta)$.

\quad

  \textbf{Theorem 2.} \emph{Suppose that
$\omega (2\delta)<2\omega (\delta)$ for all sufficiently small
$\delta>0$. There exists a nowhere linear real function
$\varphi\in C^{1, \omega}(\mathbb T)$ such that
$$
\|e^{i\lambda\varphi}\|_{A(\mathbb T)} \leq c \frac{|\lambda|}{\log
|\lambda|} \chi^{-1}\bigg(\frac{(\log |\lambda|)^2}{|\lambda|}\bigg),
\qquad \lambda\in\mathbb R, \quad |\lambda|\geq 2, \leqno (\mathrm{i})
$$
and
$$
\|e^{i\lambda\varphi}\|_{A_p(\mathbb T)} \leq c_p \bigg
(\int_1^{|\lambda|} (\chi^{-1}(1/\tau))^p ~d\tau \bigg )^{1/p},
\qquad \lambda\in\mathbb R, \quad |\lambda|\geq 2, \leqno (\mathrm{ii})
$$
for all $p, ~1<p<2$. The positive constants $c, c_p$ are
independent of $\lambda$.}

 \quad

   Put $\omega(\delta)=\delta^\alpha, ~0<\alpha<1,$ in Theorem 2.
Using Corollary 1 and the trivial estimate
$\|e^{i\lambda\varphi}\|_{A_p}\geq 1, ~1\leq p\leq 2,$ we
immediately obtain the following corollary.

\quad

   \textbf{Corollary 2.} \emph{Let $0<\alpha<1$. There exists
a nowhere linear real function $\varphi\in C^{1, \alpha}(\mathbb
T)$ such that
$$
\leqno (\mathrm{i}) \quad \|e^{i\lambda\varphi}\|_{A(\mathbb T)}=
O\big(|\lambda|^{\frac{\alpha}{1+\alpha}}(\log
|\lambda|)^{\frac{1-\alpha}{1+\alpha}}\big)
$$
and
$$
\leqno (\mathrm{ii}) \quad \|e^{i\lambda\varphi}\|_{A_p(\mathbb T)}\simeq
|\lambda|^{\frac{1}{p}-\frac{1}{1+\alpha}} \qquad \textrm{for} \quad
1<p<1+\alpha;
$$
$$
\leqno \qquad \|e^{i\lambda\varphi}\|_{A_p(\mathbb T)}\simeq 1 \qquad
\textrm{for} \quad 1+\alpha<p<2;
$$
$$
\leqno \qquad \|e^{i\lambda\varphi}\|_{A_p(\mathbb T)}=O((\log
|\lambda|)^{1/p}) \qquad \textrm{for} \quad p=1+\alpha.
$$}

\quad

  Thus, for $p\neq 1$ the estimate in Corollary 1 is
sharp. The function $\varphi\in C^{1, \alpha}$ from Corollary 2
yields the slowest possible growth of the norms
$\|e^{i\lambda\varphi}\|_{A_p}$ for $1<p<2, ~p\neq 1+\alpha$.

\quad

   \emph{Remark} 3. It is not clear if for $1<p<2$ there is
a non-trivial function $\varphi \in C^{1, p-1}(\mathbb T)$ such
that $\|e^{i\lambda\varphi}\|_{A_p(\mathbb T)}=O(1)$.

\quad

  As we mentioned in Introduction, Theorem 2 implies also the
result earlier obtained by the author in [9], namely, we have the
following assertion.

\quad

\textbf{Corollary 3.} \emph{Let $\gamma(\lambda)\geq 0$ and
$\gamma(\lambda)\rightarrow +\infty$ as $\lambda\rightarrow
+\infty$. There exists a nowhere linear real function $\varphi\in
C^1(\mathbb T)$ such that}
$$
\|e^{i\lambda\varphi}\|_{A(\mathbb T)}=O\big(\gamma(|\lambda|) \log
|\lambda|\big),
\qquad |\lambda|\rightarrow\infty, \quad \lambda\in\mathbb R.
$$

\quad

\emph{Proof.} The right-hand side in estimate (i) of Theorem  2
does not exceed $c |\lambda|\chi^{-1}(1/|\lambda|) \log |\lambda|$
and it remains to note that we can choose $\omega$ so that
$|\lambda|\chi^{-1}(1/|\lambda|)$ has an arbitrarily slow growth.
The corollary is proved.

\quad

   Before we proceed to the proof of
Theorem 2 recall that, as we mentioned in Introduction, if
$\varphi$ is a piecewise linear continuous functions then we have
$\|e^{i\lambda\varphi}\|_A=O(\log |\lambda|)$ and
$\|e^{i\lambda\varphi}\|_{A_p}=O(1)$ for $p>1$. (The logarithmic
growth in $A(\mathbb T)$ is the slowest known.) So, if we want the
norms $\|e^{i\lambda\varphi}\|_{A_p}$ to grow slowly for a $C^1$
-smooth non-trivial function $\varphi$, it seems natural to
consider the primitives of Cantor staircase functions. These
primitives resemble piecewise linear functions, although they are
$C^1$ -smooth.

   To make the proof of Theorem 2 more transparent we
start by a constructing a function $\varphi\in C^{1,
\omega}(\mathbb T), ~\varphi\neq\textrm{const}$, that satisfies
estimates (i), (ii) of the theorem, without requiring nowhere
linearity. In constructing this function we obtain certain lemmas
that we shall also use to construct a nowhere linear function
$\varphi$ with the same properties.

   Choosing (if necessarily) numbers $a, b>0$ and replacing
$\omega(\delta)$ with $a\omega(b\delta)$ (this only affects the
constants $c, c_p$ in estimates (i), (ii) in Theorem 2) we can
assume that $\omega(2\pi)=1$ and $\omega(2\delta)<2\omega(\delta)$
for $0<\delta\leq 2\pi$.

   Define positive numbers $\rho_j, ~j=0, 1, 2 \ldots,$ by
$$
\omega (\rho_j)=2^{-j}. \eqno(12)
$$

   We can assume that $\rho_0=2\pi$. Since
$\omega(2\rho_{j+1})<2\omega(\rho_{j+1})=2^{-j}=\omega(\rho_j)$ we
have $2\rho_{j+1}<\rho_j, ~j=0, 1, 2, \ldots$.

   We construct a symmetric perfect set $E\subset[0, 2\pi]$
as follows (see [8, Ch. V, \S~3], [11, Ch. XIV, \S~19]). From the
closed interval $[0, 2\pi]=[0, \rho_0]$ we remove a concentric
open interval so that two remaining closed intervals have the same
length equal to $\rho_1$. From each of the remaining closed
intervals we remove the corresponding concentric open interval so
that there remain four closed intervals of length $\rho_2$ etc. At
$j$ -th step we obtain $2^j$ closed intervals $I^j_\nu, ~\nu = 1,
2, \ldots 2^j,$ of length $\rho_j$. Continuing the process we
obtain the set of remaining points:
$$
E=\bigcap_{j=0}^\infty\bigcup_{\nu=1}^{2^j} I^j_\nu.
$$

     Let $\psi$ be the Cantor type staircase function
related to the set $E$, namely, a real continuous non-decreasing
function on $[0, 2\pi]$ that takes constant values on the
intervals complementary to $E$ in $[0, 2\pi]$ (i.e. on the
connected components of the compliment $[0, 2\pi]\setminus E$) and
increases by $2^{-j}$ on each closed interval $I_\nu^j, ~\nu=1, 2,
\ldots , 2^j,$ obtained at the $j$ -th step of construction of
$E$. Note that condition (12) implies that $\psi\in
\mathrm{Lip}_\omega([0, 2\pi])$. This can be easily verified by
nearly word for word repetition of the arguments used in [8, Ch.
V, \S~3] in the case when $\omega(\delta)=\delta^\alpha$.

   By a modification of the staircase function it is easy
to obtain a function $\psi$ on $[0, 2\pi]$ with the following
properties:

\quad

\noindent 1) $\psi$ takes constant values on the intervals
complementary to $E$ in $[0, 2\pi]$;\\
2) $\psi\in \mathrm{Lip}_\omega([0, 2\pi]) $;\\
3) $\psi(0)=\psi(2\pi)=0$;\\
4) $\int_0^{2\pi} \psi(\theta) d\theta=0$;\\
5) $\max_{t\in [0, 2\pi]}|\psi(t)|=1$;\\
6) $[0, 2\pi]=J_1\cup J_2\cup J_3$, where $J_1, J_2, J_3$ are
pairwise-disjoint intervals such that $\psi$ is monotone on each
of them.

\quad

   Define a function $\varphi$ on $\mathbb T$ by
$$
\varphi(t)=\int_0^t \psi(\theta) ~d\theta, \qquad t\in[0, 2\pi].
$$
We have $\varphi(0)=\varphi(2\pi)=\varphi'(0)=\varphi'(2\pi)=0$.
Thus, $\varphi\in C^{1, \omega}(\mathbb T)$. At the same time the
function $\varphi$ is non-constant on $\mathbb T$ and is linear on
each interval complementary to $E$ in $[0, 2\pi]$.

  We claim that $\varphi$ satisfies estimate (i) of Theorem 2.

  We put
$$
\Theta(y)=\frac{y}{\log y} \chi^{-1}\bigg(\frac{(\log y)^2}{y}\bigg),
\qquad y>1. \eqno(13)
$$

\quad

\textbf{Lemma 2.} \emph{Let $g$ be a real function on $[0, 2\pi]$
linear on an interval $\Delta\subseteq [0, 2\pi]$. Then for all
$y\geq 2$ we have
$$
\sum_{|k|\leq y} |\widefu{1_\Delta e^{ig}}(k)| \leq c \log y,
$$
where $c>0$ is independent of $\Delta, ~g$ and $y$.}

\quad

\emph{Proof.} Fix $y\ge 2$. Assume that $g(t)=at+b$ for $t\in
\Delta$, where $a, ~b \in \mathbb R$. Direct calculation yields
$$
|\widefu{1_\Delta e^{ig}}(k)|= \bigg|\frac{\sin((k-a)|\Delta|/2)}{\pi
(k-a)} \bigg|\leq\frac{1}{|k-a|}, \qquad k\in\mathbb Z, \quad k\neq a.
\eqno(14)
$$
At the same time we have
$$
|\widefu{1_\Delta e^{ig}}(k)|\leq \frac{|\Delta|}{2\pi}, \qquad
k\in\mathbb Z. \eqno(15)
$$
Thus,
$$
\sum_{|k|\leq y}|\widefu{1_\Delta e^{ig}}(k)|\leq \sum_{|k|\leq y} \min
\bigg (\frac{1}{|k-a|}, ~1\bigg ).
$$

   Consider two cases: $|a|\geq 2y$ and
$|a|<2y$. In the first case for $|k|\leq y$ we have $|k-a|\geq y$.
So,
$$
\sum_{|k|\leq y} \min \bigg (\frac{1}{|k-a|}, ~1\bigg) \leq
\sum_{|k|\leq y} \frac{1}{y} \leq 3.
$$
In the second case for $|k|\leq y$ we have $|k-a|\leq 3y$. So,
$$
\sum_{|k|\leq y} \min \bigg (\frac{1}{|k-a|}, ~1\bigg ) \leq \sum_{|k -
a|\leq 3y} \min \bigg ( \frac{1}{|k-a|}, ~1\bigg ) \leq
$$
$$
\leq \sum_{|k-a|<2} 1 +\sum_{2\leq |k-a|\leq 3y}\frac{1}{|k-a|} \leq 5 +
\int_{1\leq |x-a|\leq 3y}\frac{dx}{|x-a|}=5 + 2\log (3y).
$$
The lemma is proved.

\quad

\textbf{Lemma 3.} \emph{Let $g$ be a real function on $[0, 2\pi]$
linear on each interval complimentary to $E$ in $[0, 2\pi]$. Let
$\lambda\geq 2$. Then
$$
\sum_{|k|\leq 2\lambda} |\widefu{e^{ig}}(k)| \leq c \Theta(\lambda),
$$
where $c>0$ is independent of $\lambda$ and $g$.}

\quad

\emph{Proof.} Fix a positive integer $j$, which we shall specify
later. The set $E$ can be covered by $2^j$ closed intervals of
equal length $\rho_j$. Denote there union by $F_j$. The complement
$G_j=[0, 2\pi]\setminus F_j$ is the union of $2^j-1$
pairwise-disjoint open intervals such that on each of them $g$ is
linear. Applying Lemma 2 to each of these intervals, we obtain
$$
\sum_{|k|\leq 2\lambda} |\widefu{1_{G_j}e^{ig}}(k)| \leq c_1 2^j \log
\lambda. \eqno(16)
$$

  At the same time, $|F_j|=2^j\rho_j$, so (using the Cauchy inequality
and the Parseval identity) we have
$$
\sum_{|k|\leq 2\lambda} |\widefu{1_{F_j} e^{ig}}(k)| \leq
(4\lambda+1)^{1/2} \bigg(\sum_{|k|\leq 2\lambda}
|\widefu{1_{F_j}e^{ig}}(k)|^2 \bigg )^{1/2}
$$
$$
\leq (4\lambda+1)^{1/2} \|1_{F_j}e^{ig}\|_{L^2(\mathbb
T)}=(4\lambda+1)^{1/2}(|F_j|/(2\pi))^{1/2}
$$
$$
\leq c_2 \lambda^{1/2} (2^j\rho_j)^{1/2}. \eqno(17)
$$

    From (16), (17) we obtain
$$
\sum_{|k|\leq 2\lambda}|\widefu{e^{ig}}(k)|\leq  c_1 2^j \log \lambda
+c_2 \lambda^{1/2}(2^j\rho_j)^{1/2}. \eqno(18)
$$
The constants $c_1, c_2>0$ are independent of $\lambda, g$ and
$j$.

  Note that
$$
\log y = o(\Theta (y)), \quad y\rightarrow +\infty. \eqno(19)
$$
So if $\lambda$ is large enough, say $\lambda\geq \lambda_0$, then
$$
\frac{\Theta (\lambda)}{\log \lambda}\geq 1.
$$

  It suffices to obtain the estimate of the lemma
under additional assumption that $\lambda\geq \lambda_0$.

  Let $\lambda\geq \lambda_0$. We shall choose $j$ to
minimize the right-hand side in (18). Define a positive integer
$j(\lambda)$ by condition
$$
2^{j(\lambda)-1}\leq\frac{\Theta (\lambda)}{\log \lambda}<2^{j(\lambda)}.
\eqno(20)
$$

   It is easy to verify that
$$
\frac{\log y}{\Theta (y)}=\omega\bigg(\frac{\log y}{y}\Theta (y)\bigg),
\qquad y>1.
$$
So, the right-hand side inequality in (20) yields (see (12))
$$
\omega(\rho_{j(\lambda)})=2^{-j(\lambda)}\leq \frac{\log \lambda}{\Theta
(\lambda)}=\omega\bigg(\frac{\log \lambda}{\lambda}\Theta (\lambda)\bigg),
$$
whence
$$
\rho_{j(\lambda)}\leq \frac{\log \lambda}{\lambda}\Theta (\lambda).
\eqno(21)
$$

   The left-hand side inequality in (20) and estimate (21)
imply that for $j=j(\lambda)$ the right-hand side of (18) does not
exceed $c \Theta (\lambda)$. The lemma is proved.

\quad

   \textbf{Lemma 4.} \emph{Let $f$ be a real function on
$[0, 2\pi]$. Let $I\subseteq [0, 2\pi]$ be an interval. Suppose
that $f\in C^1(I)$, the derivative $f'$ is monotone on $I$, and
$|f'(t)|\leq 1$ for $t\in I$. Let $\lambda>0$. Then for all
$k\in\mathbb Z$ such that $|k|\geq 2\lambda$, we have}
$$
|\widefu{1_I e^{i\lambda f}}(k)|\leq \frac{2}{|k|}.
$$

\quad

   \emph{Proof.} According to the first part of the van der
Corput lemma (see [8, Ch. V, Lemma 4.3]), if $I \subset \mathbb R$
is a bounded interval and $g\in C^1(I)$ is a real function such
that its derivative $g'$ is monotone and $|g'(t)|\geq \rho>0$ for
all $t\in I$, then
$$
\bigg|\int_I e^{ig(t)} dt \bigg|\leq \frac{2\pi}{\rho}.
$$

   Fix $\lambda>0$ and $k\in\mathbb Z, ~|k|\geq 2\lambda$. Put
$g(t)=\lambda f(t)-kt$. We see that the derivative $g'(t)=\lambda
f'(t)-k$ is monotone on $I$ and $|g'(t)| \geq |k|-\lambda\geq
|k|/2$. Hence,
$$
|\widefu{1_I e^{i\lambda f}}(k)|=
\bigg|\frac{1}{2\pi}\int_I e^{ig(t)} dt \bigg|\leq \frac{2}{|k|}.
$$
The lemma is proved.

\quad

   Let us estimate the norm $\|e^{i\lambda\varphi}\|_A$. We can
assume that $\lambda\geq 2$. We shall separately estimate the sums
of the moduli of the Fourier coefficients
$\widefu{e^{i\lambda\varphi}}(k)$ over $k\in\mathbb Z$ in the
ranges $|k|\leq 2\lambda, ~2\lambda<|k|\leq\lambda^2$, and
$\lambda^2<|k|$.

   Applying Lemma 3 to the function $g=\lambda \varphi$, we obtain
$$
\sum_{|k|\leq 2\lambda} |\widefu{e^{i\lambda \varphi}}(k)| \leq
c_1\Theta (\lambda). \eqno(22)
$$

  Note that $|\varphi'(t)|=|\psi(t)|\leq 1$ on
$[0, 2\pi]$ and the interval $[0, 2\pi]$ is a union of three
pairwise-disjoint intervals such that the derivative
$\varphi'=\psi$ is monotone on each of them. Applying Lemma 4 to
each of these intervals, we see that
$$
|\widefu{e^{i\lambda\varphi}}(k)|\leq \frac{6}{|k|}, \qquad
|k|>2\lambda. \eqno(23)
$$
So,
$$
\sum_{2\lambda<|k|\leq\lambda^2} |\widefu{e^{i\lambda \varphi}}(k)|\leq
\sum_{2\lambda<|k|\leq\lambda^2} \frac{6}{|k|}\leq c_2 \log \lambda.
\eqno(24)
$$

   Note then that
$$
\sum_{|k|>\lambda^2} |\widefu{e^{i\lambda
\varphi}}(k)|=\sum_{|k|>\lambda^2}\frac{1}{|k|} |\widefu{(e^{i\lambda
\varphi})'}(k)|\leq\bigg(\sum_{|k|>\lambda^2}
\frac{1}{k^2}\bigg)^{1/2} \|(e^{i\lambda
\varphi})'\|_{L^2(\mathbb T)}\leq c_3. \eqno(25)
$$

  Summing inequalities (22), (24), (25)
and taking (19) into account, we obtain estimate (i) of Theorem 2.

\quad

   We claim that the same function
$\varphi$ satisfies estimate (ii) of Theorem 2. For $p>1$ we put
$$
\Theta_p(y)=\bigg (\int_1^{y} (\chi^{-1}(1/\tau))^p ~d\tau \bigg )^{1/p},
\qquad y>1. \eqno(26)
$$

\quad

   \textbf{Lemma 5.} \emph{Let $1<p<2$. Let $g$ be a real function
on $[0, 2\pi]$ linear on an interval $\Delta\subseteq [0, 2\pi]$.
Then
$$
\|1_\Delta e^{ig}\|_{A_p(\mathbb T)}\leq c_p |\Delta|^{1/q},
$$
where $1/p+1/q=1$ and the constant $c_p>0$ is independent of $g$
and $\Delta$.}

\quad

   \emph{Proof.} Assuming that $g(t)=at+b$ for
$t\in \Delta$, from estimates (14), (15), used in the proof of
Lemma 2, we obtain
\begin{multline*}
\|1_\Delta e^{ig}\|_{A_p(\mathbb T)}^p=\sum_{k\in\mathbb
Z}|\widefu{1_\Delta e^{ig}}(k)|^p \leq \sum_{|k-a|>1/|\Delta|}
\frac{1}{|k-a|^p}+\sum_{|k-a|\leq 1/|\Delta|} |\Delta|^p\leq c_p
|\Delta|^{p-1}.
\end{multline*}
The lemma is proved.

\quad

   In what follows it will be convenient to use the
analogues of the spaces $A_p(\mathbb T)$ for the functions defined
on the line $\mathbb R$. We use the same symbol ${}^\wedge$ to
denote the Fourier transform of tempered distributions on $\mathbb
R$. For $1<p<\infty$ let $A_p(\mathbb R)$ be the space of tempered
distributions  $g$ on $\mathbb R$ such that $\fu{g}$ belongs to
$L_p(\mathbb R)$. We put
$$
\|g\|_{A_p(\mathbb R)}=\|\fu{g}\|_{L^p(\mathbb R)}.
$$

   It is known (see, e.g. [12]) that for $1\leq p\leq 2$
the Fourier transform (as well as its inverse) is a bounded
operator from $L^p(\mathbb R)$ to $L^q(\mathbb R), ~1/p+1/q =1$.
Thus, each distribution that belongs to $A_p(\mathbb R), ~1<p\leq
2,$ is a function in $L^q(\mathbb R)$.

  Let $1<p\leq 2$. It is easy to verify (and is well known,
see, e.g. [13, \S~44]) that if $f$ is a $2\pi$ -periodic function
and $f^\ast$ is its restriction to $[0, 2\pi]$ extended by zero to
$\mathbb R$, that is $f^\ast=f$ on $[0, 2\pi], ~f^\ast=0$ on
$\mathbb R\setminus [0, 2\pi]$, then $f$ belongs to $A_p(\mathbb
T)$ if and only if  $f^\ast\in A_p(\mathbb R)$. The norms satisfy
$$
c_1(p) \|f^\ast\|_{A_p(\mathbb R)}\leq \|f\|_{A_p(\mathbb T)}\leq c_2(p)
\|f^\ast\|_{A_p(\mathbb R)}.
$$
We call this statement the transference principle.

   Let $\Delta\subseteq\mathbb R$ be an arbitrary interval.
It is known (see, e.g., [14]) that the operator $S_\Delta$ given
by
$$
S_\Delta (g)=(1_\Delta \cdot \fu{g})^\rightuf,
\qquad g\in L^p\cap L^2(\mathbb R),
$$
where ${}^\vee$ means the inverse Fourier transform, is a bounded
operator from $L^p(\mathbb R)$ to itself for $1<p<\infty$.

   For an arbitrary family $\{\Delta_\nu\}$ of pairwise-disjoint
intervals in $\mathbb R$ we define the Littlewood-Paley square
function:
$$
S(g)=\bigg (\sum_\nu |S_{\Delta_\nu} (g)|^2\bigg)^{1/2}.
$$

  We recall the Rubio de Francia inequality [15]:
$$
\|S(g)\|_{L^p(\mathbb R)}\leq c_p \|g\|_{L^p(\mathbb R)}, \qquad g\in
L^p(\mathbb R),
$$
which holds for $2<p<\infty$. The constant $c_p>0$ is independent
of $g$ and $\{\Delta_\nu\}$. By duality, for $1<p<2$ and for any
function $g\in L^p(\mathbb R)$ with
$$
\textrm{supp}~\fu{g}~\subseteq \bigcup_\nu \Delta_\nu, \eqno(27)
$$
we have
$$
\|g\|_{L^p(\mathbb R)}\leq c_p \|S(g)\|_{L^p(\mathbb R)}. \eqno(28)
$$

   The following lemma is a simple consequence
of the Rubio de Francia inequality.

\quad

   \textbf{Lemma 6.} \emph{Let $1<p<2$. Let
$\Delta_\nu\subset\mathbb T, ~\nu=1, 2, \ldots, N,$ be a finite
family of pairwise-disjoint intervals and let $f_\nu, ~\nu=1, 2,
\ldots, N,$ be functions in $A_p(\mathbb T)$ such that
$\mathrm{supp} \,f_\nu\subseteq\Delta_\nu, ~\nu=1, 2, \ldots , N$.
Then
$$
\bigg \|\sum_\nu f_\nu \bigg \|_{A_p(\mathbb T)}^p\leq c_p\sum_\nu
\|f_\nu\|_{A_p(\mathbb T)}^p,
$$
where $c_p>0$ is independent of $N$ and the families
$\{\Delta_\nu\}, ~\{f_\nu\}$.}

\quad

   \emph{Proof.} Consider a finite family
$\Delta_\nu, ~\nu=1, 2, \ldots, N,$ of pairwise-disjoint bounded
intervals in $\mathbb R$. For $1<p<2$, taking into account the
obvious inequality
$$
S(g)\leq \bigg (\sum_\nu |S_{\Delta_\nu} (g)|^p\bigg)^{1/p},
$$
we obtain from (28) that
$$
\|g\|_{L^p(\mathbb R)}^p\leq c_p^p \sum_\nu
\|S_{\Delta_\nu}(g)\|_{L^p(\mathbb R)}^p \eqno(29)
$$
for any function $g\in L^p(\mathbb R)$ with condition (27).

  Let $f_\nu, ~\nu=1, 2, \ldots, N,$ be functions that
belong to $A_p(\mathbb R), ~1<p< 2,$ and have supports in
$\Delta_\nu, ~\nu=1, 2, \ldots, N,$ respectively. Applying
estimate (29) to the function
$$
g=\bigg(\sum_\nu f_\nu \bigg)^\rightuf,
$$
and taking into account that the direct and inverse Fourier
transforms on $\mathbb R$ differ only in the sign of the variable,
we see that
$$
\bigg \|\sum_\nu f_\nu \bigg \|_{A_p(\mathbb R)}^p\leq c_p^p\sum_\nu
\|f_\nu\|_{A_p(\mathbb R)}^p.
$$
(Since the functions $f_\nu$ are in $L^q(\mathbb R)$ and the
intervals $\Delta_\nu$ are bounded, we have $f_\nu\in L^2(\mathbb
R), ~\nu=1, 2, \ldots, N$. So, $g\in L^2(\mathbb R)$.)

   It remains to use the transference principle.
The lemma is proved.

\quad

   \textbf{Lemma 7.} \emph{Let $1<p<2$ and let $g$ be a real function
on $[0, 2\pi]$ linear on each interval complementary to $E$ in
$[0, 2\pi]$. Let $\lambda\geq 2$. Then
$$
\sum_{|k|\leq 2\lambda} |\widefu{e^{ig}}(k)|^p \leq c_p
\,(\Theta_p(\lambda))^p,
$$
where $c_p>0$ is independent of $\lambda$ and $g$.}

\quad

  \emph{Proof.} For an arbitrary $j=1, 2, \ldots$ the set $E$
can be covered by $2^j$ closed intervals of equal length $\rho_j$.
Their union $F_j$ is of measure $|F_j|=2^j\rho_j$. The compliment
$G_j=[0, 2\pi]\setminus F_j$ is a union of $2^j-1$
pairwise-disjoint open intervals. Denote these intervals by
$\Delta_\nu^j, ~\nu=1, 2, \ldots 2^j-1$. The function $g$ is
linear on each interval $\Delta_\nu^j$, so, using Lemmas 6 and 5,
we obtain that
$$
\|1_{G_j}e^{ig}\|_{A_p(\mathbb T)}^p\leq c_p \sum_{\nu=1}^{2^j-1}
\|1_{\Delta_\nu^j}e^{ig}\|_{A_p(\mathbb T)}^p \leq c_{p, 1}
\sum_{\nu=1}^{2^j-1}|\Delta_\nu^j|^{p-1}. \eqno(30)
$$

  The family of intervals $\Delta_\nu^j, ~\nu=1,
2, \ldots, 2^j-1,$ consists of one interval of length
$\rho_0-2\rho_1<\rho_0$, of two intervals of length
$\rho_1-2\rho_2<\rho_1$, etc., of $2^m$ intervals of length
$\rho_m-2\rho_{m+1}<\rho_m$, where $m=0, 1, \ldots j-1$. Hence we
obtain (see (30)) that
$$
\sum_{|k|\leq 2\lambda} |\widefu{1_{G_j}e^{ig}}(k)|^p \leq
\|1_{G_j}e^{ig}\|_{A_p(\mathbb T)}^p\leq c_{p, 1} \sum_{m=0}^{j-1}2^m
\rho_m^{p-1}. \eqno(31)
$$

  At the same time, using the H\"older inequality with
$p^*=2/p$ and $1/p^*+1/q^*=1$, we have
$$
\sum_{|k|\leq 2\lambda} |\widefu{1_{F_j}e^{ig}}(k)|^p \leq \bigg
(\sum_{|k|\leq 2\lambda} |\widefu{1_{F_j}e^{ig}}(k)|^{pp^*} \bigg
)^{1/p^*} \bigg (\sum_{|k|\leq 2\lambda}1 \bigg )^{1/q^*}
$$
$$
\leq\bigg (\sum_{|k|\leq 2\lambda} |\widefu{1_{F_j}e^{ig}}(k)|^2 \bigg
)^{p/2} (4\lambda+1)^{1-\frac{p}{2}}\leq
\|1_{F_j}e^{ig}\|^p_{L^2(\mathbb
T)}(4\lambda+1)^{1-\frac{p}{2}}
$$
$$
=(|F_j|/(2\pi))^{p/2}(4\lambda+1)^{1-\frac{p}{2}}\leq c_{p, 2}(2^j
\rho_j)^{p/2}\lambda^{1-\frac{p}{2}}. \eqno(32)
$$

   From (31), (32) we obtain that
$$
\bigg (\sum_{|k|\leq 2\lambda} |\widefu{e^{ig}}(k)|^p \bigg )^{1/p}\leq
c_{p, 3} \bigg (\sum_{m=0}^{j-1}2^m \rho_m^{p-1} \bigg )^{1/p}+c_{p,
4}(2^j \rho_j)^{1/2}\lambda^{\frac{1}{p}-\frac{1}{2}}, \eqno(33)
$$
where $c_{p, 3}, c_{p, 4}>0$ are independent of $\lambda, g$ and
$j$.

   It is easy to choose $j$ so that to minimize the expression
on the right-hand side of (33) in the case when $\omega
(\delta)=\delta^\alpha$. General case requires certain
calculations.

   We put
$$
a_m=\frac{1}{\chi (\rho_m)}, \qquad m=0, 1, 2, \ldots.
$$
The sequence $\{a_m\}$ is unbounded, strictly increases, and
$a_0=1/(2\pi)$.

   It is clear that it suffices to obtain the
estimate of the lemma under additional condition $\lambda\geq
a_1$.

   Choose $j(\lambda)=2, 3, \ldots$ so that
$$
a_{j(\lambda)-1}\leq \lambda<a_{j(\lambda)}. \eqno(34)
$$

   Let us estimate the first term on the right-hand side
in (33). Note that for $m=1, 2, \ldots$
$$
\frac{a_{m}}{a_{m-1}}=\frac{\chi (\rho_{m-1})}{\chi
(\rho_m)}=\frac{\rho_{m-1}\omega(\rho_{m-1})}{\rho_m\omega(\rho_m)}=
\frac{\rho_{m-1}2^{-(m-1)}}{\rho_m 2^{-m}}=\frac{2\rho_{m-1}}{\rho_m}\geq
2,
$$
whence
$$
a_m\leq 2(a_m-a_{m-1}), \qquad m=1, 2, \ldots.
$$
So,
$$
\sum_{m=1}^{j(\lambda)-1}2^m
\rho_m^{p-1}=\sum_{m=1}^{j(\lambda)-1}
\frac{\rho_m^{p-1}}{\omega(\rho_m)}=
\sum_{m=1}^{j(\lambda)-1}\frac{\rho_m^p}{\chi(\rho_m)}
$$
$$
=\sum_{m=1}^{j(\lambda)-1}a_m(\chi^{-1}(1/a_m))^p
\leq 2\sum_{m=1}^{j(\lambda)-1}(a_m-a_{m-1})(\chi^{-1}(1/a_m))^p
$$
$$
\leq 2\sum_{m=1}^{j(\lambda)-1}\int_{a_{m-1}}^{a_m}
(\chi^{-1}(1/\tau))^p~d\tau
=2\int_{a_0}^{a_{j(\lambda)-1}}(\chi^{-1}(1/\tau))^p~d\tau
$$
$$
\leq 2\int_{a_0}^\lambda(\chi^{-1}(1/\tau))^p~d\tau=c_{p,
5}+2(\Theta_p(\lambda))^p. \eqno(35)
$$

  Let us estimate the second term on the right-hand side
in (33). Put $t=\chi^{-1}(1/\lambda)$. Then (see (34))
$$
\rho_{j(\lambda)}<t\leq \rho_{j(\lambda)-1}.
$$
Hence, using the right-hand side inequality, we obtain that
$\omega(t)\leq 2^{-(j(\lambda)-1)}$ and therefore
$$
2^{j(\lambda)}\leq
\frac{2}{\omega(t)}=\frac{2t}{\chi(t)}=2\chi^{-1}(1/\lambda)\lambda.
$$
At the same time, the left-hand side inequality yields
$$
\rho_{j(\lambda)}<t=\chi^{-1}(1/\lambda).
$$
So,
$$
(2^{j(\lambda)}\rho_{j(\lambda)})^{p/2}\lambda^{1-\frac{p}{2}}\leq
2^{p/2}\lambda (\chi^{-1}(1/\lambda))^p.
$$
Since $\lambda\geq 2$ and the function $\chi^{-1}$ is increasing,
we obtain
$$
\lambda(\chi^{-1}(1/\lambda))^p\leq
2(\lambda-1)(\chi^{-1}(1/\lambda))^p\leq 2\int_1^\lambda
(\chi^{-1}(1/\tau))^p d\tau=2(\Theta_p(\lambda))^p.
$$
Thus
$$
(2^{j(\lambda)}\rho_{j(\lambda)})^{p/2}\lambda^{1-\frac{p}{2}}\leq c_{p,
6}(\Theta_p(\lambda))^p.
$$

   Using this inequality and (35), we see that
for $j=j(\lambda)$ the right-hand side in (33) does not exceed
$c_p \Theta_p(\lambda)$. The lemma is proved.

\quad

   Let us estimate the norms $\|e^{i\lambda\varphi}\|_{A_p}, ~p>1$.
Applying Lemma 7 to the function $g=\lambda\varphi$, we have (we
can assume that $\lambda\ge 2$)
$$
\sum_{|k|\leq 2\lambda} |\widefu{e^{i\lambda \varphi}}(k)|^p \leq c_{p,
1} (\Theta_p (\lambda))^p.
$$
   At the same time, using estimate (23), we have
$$
\sum_{|k|>2\lambda} |\widefu{e^{i\lambda \varphi}}(k)|^p \leq
\sum_{|k|>2\lambda}\bigg(\frac{6}{|k|}\bigg)^p \leq 6^p \sum_{|k|\geq 1}
\frac{1}{|k|^p} = c_{p, 2}.
$$
Thus,
$$
\|e^{i\lambda \varphi}\|^p_{A_p(\mathbb T)}\leq c_{p, 1}(\Theta_p
(\lambda))^p+c_{p, 2},
$$
and we obtain (ii).

\quad

  \emph{Proof of Theorem} 2. For $m=0, 1, 2, \ldots$
let $E_m$ be the portion of the set $E$ in the closed interval
$[0, \rho_m]$, i.e., $E_m=E\cap [0, \rho_m]$. Certainly $E_0=E$.

   Again let $\psi$ be the Cantor
type staircase function related to the set $E$. For each $m=0, 1,
2, \ldots$ the restriction of $\psi$ to $[0, \rho_m]$ belongs to
$\mathrm{Lip}_\omega([0, \rho_m])$. Thus, it is easy to see that
there exist functions $\psi_m, ~m=0, 1, 2, \ldots,$ on $[0, 2\pi]$
(modified staircase functions) with the following properties:

\quad

\noindent 1) $\psi_m$ takes constant values on the intervals
complementary to $E_m$ in $[0, \rho_m]$;\\
2) $\psi_m\in \mathrm{Lip}_\omega([0, \rho_m]) $;\\
3) $\psi_m(0)=\psi_m(\rho_m)=0$ and $\psi_m=0$ on
$[0, 2\pi]\setminus [0, \rho_m]$;\\
4) $\int_0^{\rho_m} \psi_m (\theta) d\theta=0$;\\
5) $\max_{t\in [0, \rho_m]}|\psi_m(t)|=1$;\\
6) $[0, \rho_m]$ is a union of three pairwise-disjoint intervals
such that on each of them $\psi_m$ is monotone.

\quad

   For $m=0, 1, 2, \ldots$ we put
$$
\varphi_m(t)=\int_0^t \psi_m (\theta)~d\theta, \qquad t\in [0, 2\pi].
$$
It is clear that
$\varphi_m(0)=\varphi_m(\rho_m)=\varphi_m'(0)=\varphi_m'(\rho_m)=0$,
$\varphi_m=0$ on $[0, 2\pi]\setminus [0, \rho_m]$ and
$\varphi_m'\in \mathrm{Lip}_\omega [0, \rho_m], ~m=0, 1, 2,
\ldots$.

     Let $I\subseteq [0, 2\pi]$ be a closed interval.
Let $E_m(I)$ be an affine copy of the set $E_m$ related to $I$,
namely: $E_m(I)=l_{m, I}^{-1}(E_m)$, where $l_{m, I}$ is an affine
map that maps $I$ onto $[0, \rho_m]$. Let $\varphi_m^I$ be an
affine copy of the function $\varphi_m$ related to $I$, namely a
function on $[0, 2\pi]$ such that $\varphi_m^I=\varphi_m\circ
l_{m, I}$ on $I$ and $\varphi_m^I=0$ on $[0, 2\pi]\setminus I$.

    Using induction we define a sequence of
closed intervals $I_m\subseteq [0, 2\pi]$ with the condition
$$
|I_m|\leq 2^{-m}\rho_m, \qquad m=0, 1, 2, \ldots, \eqno(36)
$$
and a sequence of sets $B_m\subset I_m, ~m=0, 1, 2, \ldots$. Let
$I_0=[0, 2\pi]$ and $B_0=E_0$. Once $I_0, I_1, I_2, \ldots , I_m$
and $B_0, B_1, B_2, \ldots, B_m$ have been defined, we define
$I_{m+1}$ and $B_{m+1}$ as follows. Consider the union
$\bigcup_{s=1}^m B_s$ and choose an interval of maximum length
complementary to this union in $[0, 2\pi]$ (if there are several
of them, we take any one). Denote this interval by $J$. Let
$I_{m+1}$ be a closed interval contained in $J$, concentric with
$J$, and of length $|I_{m+1}|\leq 2^{-(m+1)}\rho_{m+1}$. We put
$B_{m+1}=E_{m+1}(I_{m+1})$.

   Now we define functions $f_m, ~m=0, 1, 2, \ldots,$ by
$f_m=\varphi_m^{I_m}$.

   In what follows we consider $2\pi$ -periodic extensions of
functions $\varphi_m, ~f_m, ~m=0, 1, 2, \ldots,$ from $[0, 2\pi]$
to $\mathbb R$ and thus we regard these functions as functions on
$\mathbb T$. All these functions belong to $C^{1, \omega}(\mathbb
T)$. Each function $f_m$ vanishes on $[0, 2\pi]\setminus I_m$ and
is linear on the intervals complimentary to $B_m$ in $I_m$. We
also have
$$
|\varphi_m'(t)|\leq 1, \qquad t\in\mathbb T, \quad m=0, 1, 2, \ldots,
$$
and, since $l_{m, I_m}$ is a linear function on $[0, 2\pi]$ with
tangent coefficient equal (up to a sign) to $\rho_m/|I_m|$, we see
that
$$
|f_m'(t)|\leq \rho_m/|I_m|, \qquad t\in\mathbb T, \quad m=0, 1, 2,
\ldots. \eqno(37)
$$

   We put
$$
\varphi=\sum_{m=0}^\infty \varepsilon_m f_m,
$$
where numbers $\varepsilon_m>0$ decrease to $0$ so fast that
$\varphi\in C^{1, \omega}(\mathbb T)$.

   It is clear that the function $\varphi$ is nowhere linear.

   In addition we require that the sequence
$\{\varepsilon_m\}$ satisfies
$$
\sum_{m=0}^\infty \varepsilon_m \frac{\rho_m}{|I_m|}=1 \eqno(38)
$$
and decreases so fast that
$$
\delta_{j(\lambda)}=O(1/\lambda),
\quad \lambda\rightarrow +\infty, \eqno(39)
$$
where
$$
\delta_j=\sum_{m=j+1}^\infty \varepsilon_m \frac{\rho_m}{|I_m|} \eqno(40)
$$
and $j(\lambda)$ is the positive integer defined for each
sufficiently large $\lambda$ by condition (20) (it is clear that
$j(\lambda)\rightarrow\infty$ as $\lambda\rightarrow +\infty$, see
(19)).

  Let us estimate the norms $\|e^{i\lambda\varphi}\|_A$.
We can assume that $\lambda\geq 2$.

  Put
$$
S_j=\sum_{m=0}^j \varepsilon_m f_m, \qquad j=0, 1, 2, \ldots.
$$

 For each $j$ the set $E_m, ~m\leq j,$ can be covered by
$2^{j-m}$ closed intervals of length $\rho_j$. Since $B_m$ is
obtained from $E_m$ by contraction (see (36)), we see that the set
$B_m, ~m\leq j,$ can be covered by $2^{j-m}$ closed intervals
whose length does not exceed $\rho_j$. Therefore, the set
$$
\bigcup_{m=0}^j B_m \eqno(41)
$$
can be covered by
$$
2^j+2^{j-1}+2^{j-2}+\ldots +2^1+1\leq 2^{j+1}
$$
closed intervals of length $\rho_j$. Denote the union of these
closed intervals by $F_j$. We have
$$
|F_j|\leq 2^{j+1}\rho_j.
$$

   We put $G_j=[0, 2\pi]\setminus F_j$. It is clear that $G_j$
is a union of at most $2^{j+1}+1$ pairwise-disjoint intervals. It
is clear that the function $S_j$ is linear on each interval
complementary to the union (41) in $[0, 2\pi]$. Therefore the
function $\lambda S_j$ is linear on the intervals that form $G_j$.
Applying Lemma 2 to each of these intervals, we have
$$
\sum_{|k|\leq 2\lambda} |\widefu{1_{G_j}e^{\lambda S_j}}(k)|\leq c_1 2^j
\log\lambda. \eqno(42)
$$

  At the same time,
$$
\sum_{|k|\leq 2\lambda} |\widefu{1_{F_j}e^{\lambda S_j}}(k)|\leq
(4\lambda+1)^{1/2}\bigg (\sum_{|k|\leq 2\lambda}
|\widefu{1_{F_j}e^{\lambda S_j}}(k)|^2\bigg )^{1/2}
$$
$$
\leq (4\lambda+1)^{1/2}\|1_{F_j}e^{\lambda S_j}\|_{L^2(\mathbb T)}=
(4\lambda+1)^{1/2}(|F_j|/(2\pi))^{1/2}
$$
$$
\leq c_2 \lambda^{1/2} (2^j\rho_j)^{1/2}. \eqno(43)
$$

  By summing estimates (42) and (43) we obtain
$$
\sum_{|k|\leq 2\lambda} |\widefu{e^{\lambda S_j}}(k)|\leq c_1 2^j
\log\lambda+c_2\lambda^{1/2} (2^j\rho_j)^{1/2}, \eqno(44)
$$
where $c_1, c_2>0$ are independent of $\lambda$ and $j$.

  We note that (see (37), (38)),
$$
|S_j'(t)|\leq \sum_{m=0}^j \varepsilon_m \frac{\rho_m}{|I_m|}\leq 1,
\qquad t\in \mathbb T. \eqno(45)
$$

  We also note that for every $j=0, 1, 2, \ldots$ the
interval $[0, 2\pi]$ is a union of three pairwise disjoint
intervals such that the derivative $f_j'$ of $f_j$ is monotone on
each of them. In addition, $f_j'$ is supported in a certain
interval on which the derivative $f_m'$ of any function $f_m,
~m<j,$ is constant. So, for every $j=0, 1, 2, \ldots$ the interval
$[0, 2\pi]$ is a union of at most $4j+3$ pairwise disjoint
intervals such that on each of them the derivative of $S_j$ is
monotone (this can be easily verified by induction with respect to
$j$).

  Taking into account estimate (45) and applying Lemma 4 to
each interval on which $S_j'$ is monotone, we obtain for
$|k|>2\lambda$
$$
|\widefu{e^{i\lambda S_j}}(k)|\leq \frac{8j+6}{|k|},
$$
whence
$$
\sum_{2\lambda<|k|\leq \lambda^2} |\widefu{e^{\lambda S_j}}(k)|\leq c_3
\,(8j+6)\log\lambda. \eqno(46)
$$

  Note then, that we have $|S_j'(t)|\leq 1$ for $t\in\mathbb T$ (see (45)),
so,
$$
\sum_{|k|>\lambda^2} |\widefu{e^{\lambda S_j}}(k)|=\sum_{|k|>\lambda^2}
\frac{1}{|k|}|\widefu{(e^{\lambda S_j})'}(k)|\leq
\bigg(\sum_{|k|> \lambda^2}
\frac{1}{|k|^2}\bigg)^{1/2}\|(e^{\lambda S_j})'\|_{L^2(\mathbb
T)}\leq c_4. \eqno(47)
$$

  From (44), (46), (47) we obtain
$$
\|e^{i\lambda S_j}\|_{A(\mathbb T)}\leq c_5 2^j
\log\lambda+c_6\lambda^{1/2} (2^j\rho_j)^{1/2}, \eqno(48)
$$
where $c_5, c_6>0$ are independent of $\lambda, j$.

  We put now
$$
r_j=\sum_{m=j+1}^\infty \varepsilon_m f_m.
$$

   Obviously, for any function $g\in C^1(\mathbb T)$ we have
$$
\|g\|_{A(\mathbb T)}\leq c\|g\|_{C^1(\mathbb T)}, \eqno(49)
$$
where
$$
\|g\|_{C^1(\mathbb T)}=\max_{t\in\mathbb T}|g(t)|+\max_{t\in\mathbb
T}|g'(t)|,
$$
and $c>0$ is independent of $g$. Since
$$
|r_j'(t)|\leq \delta_j, \qquad t\in \mathbb T
$$
(see (37), (40)), we obtain
$$
\|e^{i\lambda r_j}\|_{A(\mathbb T)}
\leq c \|e^{i\lambda r_j}\|_{C^1(\mathbb
T)}\leq c(1+\lambda \delta_j). \eqno(50)
$$

  We note now that the expression on the
right-hand side in (48) is of the same form as in (18). We
minimize it in the same way as in the proof of Lemma 3, namely,
assuming that $\lambda$ is sufficiently large we choose
$j=j(\lambda)$ satisfying condition (20). We see that for this $j$
the right-hand side in (48) does not exceed $c_7 \Theta
(\lambda)$. Thus,
$$
\|e^{i\lambda S_{j(\lambda)}}\|_{A(\mathbb T)}\leq c_7 \Theta (\lambda).
\eqno(51)
$$

   At the same time, from (50), taking into account
condition (39), we obtain
$$
\|e^{i\lambda r_{j(\lambda)}}\|_{A(\mathbb T)}\leq c_8. \eqno(52)
$$

   Since $\varphi=S_j+r_j$ for any $j$, we have
$$
\|e^{i\lambda \varphi}\|_{A(\mathbb T)} \leq \|e^{i\lambda
S_j}\|_{A(\mathbb T)} \|e^{i\lambda r_j}\|_{A(\mathbb T)}.
$$
Thus, from (51), (52) we obtain estimate (i) of Theorem 2.

\quad

  Let us estimate the norms $\|e^{i\lambda\varphi}\|_{A_p}, ~p>1$.
We need two simple lemmas.

\quad

\textbf{Lemma 8.} \emph{Let $1<p<2$. Let $I, J$ be two intervals
in $[0, 2\pi]$ and let $U, V$ be two functions in $A_p(\mathbb T)$
vanishing on $[0, 2\pi]\setminus I$ and $[0, 2\pi]\setminus J$
respectively. Suppose that $U(t)=V\circ l(t)$ for $t\in I$, where
$l$ is an affine map of $I$ onto $J$. Then $\|U\|_{A_p(\mathbb
T)}\leq c_p|a|^{-1/q}\|V\|_{A_p(\mathbb T)}$, where $a$ is the
tangent coefficient of  $l, ~|a|=|J|/|I|,$ and $1/p+1/q=1$. The
constant $c_p>0$ is independent of $U, V, l, I, J$.}

\quad

\emph{Proof.} If $l(t)=at+b, ~t\in\mathbb R, ~a\neq 0,$ then for
an arbitrary function $g\in L^1(\mathbb R)\cap A_p(\mathbb R)$
direct calculation yields
$$
|\widefu{g\circ l}(u)|=
\bigg|\frac{1}{a}\fu{g}\bigg(\frac{u}{a}\bigg)\bigg|,
$$
and therefore,
$$
\|g\circ l\|_{A_p(\mathbb R)}=\frac{1}{|a|^{1/q}}\|g\|_{A_p(\mathbb R)}.
$$
It remains to use the transference principle. The lemma is proved.

\quad

\textbf{Lemma 9.} \emph{Let $1<p<2$ and let $l$ be a real $2\pi$
-periodic function linear on $(0, 2\pi)$. Then for any function
$f\in A_p(\mathbb T)$ we have $e^{il} f\in A_p(\mathbb T)$ and
$\|e^{il}f\|_{A_p(\mathbb T)}\leq c_p\|f\|_{A_p(\mathbb T)}$ where
$c_p>0$ is independent of $f$ and $l$.}

\quad

\emph{Proof.} For an arbitrary linear function  $l(t)=at+b$ on
$\mathbb R$ and for an arbitrary function $g\in L^1(\mathbb R)\cap
A_p(\mathbb R)$ we have
$$
|\widefu{e^{il} g}(u)|=|\fu{g}(u-a)|, \qquad u\in\mathbb R,
$$
whence $e^{il}g\in A_p(\mathbb R)$ and $\|e^{il} g\|_{A_p(\mathbb
R)}=\|g\|_{A_p(\mathbb R)}$. It remains to use the transference
principle. The lemma is proved.

\quad

  Let $\lambda\geq 2$. For each $m=0, 1, 2 \ldots$
the function $\varphi_m$ is linear on the intervals complementary
to $E_m$ in $[0, 2\pi]$ and therefore (since $E_m\subseteq
E_0=E$), on the intervals complementary to $E$ in $[0, 2\pi]$.
Using Lemma 7, we obtain
$$
\sum_{|k|\leq 2\lambda} |\widefu{e^{i\lambda \varphi_m}}(k)|^p\leq c_{p,
1} (\Theta_p(\lambda))^p, \qquad \lambda \geq 2, \quad m=0, 1, 2 \ldots.
\eqno(53)
$$

  For each $m=0, 1, 2, \ldots$ we have
$|\varphi_m'(t)|\leq 1, ~t\in [0, 2\pi],$ and the interval $[0,
2\pi]$ can be partitioned into three intervals so that the
derivative $\varphi_m'$ is monotone on each of them. Using Lemma 4
we obtain
$$
|\widefu{e^{i\lambda \varphi_m}}(k)|\leq \frac{6}{|k|}
$$
for $|k|>2\lambda$. So,
$$
\sum_{|k|>2\lambda} |\widefu{e^{i\lambda \varphi_m}}(k)|^p\leq c_{p,2},
\qquad \lambda \geq 2, \quad m=0, 1, 2 \ldots.
$$
Together with (53) this yields
$$
\|e^{i\lambda\varphi_m}\|_{A_p(\mathbb T)}\leq c_{p, 3}\Theta_p(\lambda),
\qquad \lambda\geq 2, \quad m=0, 1, 2 \ldots. \eqno(54)
$$
Let now $0\leq\lambda<2$. We obtain that (see (49))
$$
\|e^{i\lambda \varphi_m}\|_{A_p(\mathbb T)}\leq \|e^{i\lambda
\varphi_m}\|_{A(\mathbb T)}\leq c \|e^{i\lambda \varphi_m}\|_{C^1(\mathbb
T)}\leq c(1+\lambda) \leq 3c,
$$
$$
\qquad 0\leq\lambda<2, \quad m=0, 1, 2 \ldots. \eqno(55)
$$

  In what follows we assume that $\lambda\geq 2$. It is clear that
$0<\varepsilon_m\leq 1$ for all $m$ (see (36), (38)). Applying
estimates (54), (55) in the cases when $\lambda\varepsilon_m\geq
2$ and $\lambda\varepsilon_m< 2$ respectively, we have
$$
\|e^{i\lambda \varepsilon_m \varphi_m}\|_{A_p(\mathbb T)}\leq c_{p,
4}\Theta_p(\lambda), \qquad \lambda\geq 2, \quad m=0, 1, 2, \ldots,
\eqno(56)
$$
whence
$$
\|e^{i\lambda \varepsilon_m \varphi_m}-1\|_{A_p(\mathbb T)}\leq c_{p,
5}\Theta_p(\lambda), \qquad \lambda\geq 2,
\quad m=0, 1, 2, \ldots. \eqno(57)
$$

  For a fixed $m$ we apply Lemma 8 to the closed intervals
$I=I_m, ~J=[0, \rho_m]$ and the functions $U=
e^{i\lambda\varepsilon_m f_m}-1,
~V=e^{i\lambda\varepsilon_m\varphi_m}-1$. For the corresponding
tangent coefficient $a$ we have $|a|=\rho_m/|I_m|\geq 2^m$ (see
(36)), so, from (57) we obtain that
$$
\|e^{i\lambda \varepsilon_m f_m}-1\|_{A_p(\mathbb T)}\leq c_{p,
6}2^{-m/q}\|e^{i\lambda \varepsilon_m \varphi_m}-1\|_{A_p(\mathbb T)}\leq
$$
$$
\leq c_{p, 7}2^{-m/q}\Theta_p(\lambda), \qquad \lambda\geq 2, \quad m=0,
1, 2, \ldots. \eqno(58)
$$

  We note now that for $j=0, 1, 2, \ldots$
$$
e^{i\lambda S_{j+1}}-e^{i\lambda S_j}=1_{I_{j+1}}(e^{i\lambda
S_{j+1}}-e^{i\lambda S_j})=1_{I_{j+1}}e^{i\lambda S_j}(e^{i\lambda
\varepsilon_{j+1}f_{j+1}}-1),
$$
and the function $S_j$ is linear on $I_{j+1}$. So, using Lemma 9,
we have
$$
\|e^{i\lambda S_{j+1}}-e^{i\lambda S_j}\|_{A_p(\mathbb T)}\leq c_{p,
8}\|1_{I_{j+1}}(e^{i\lambda \varepsilon_{j+1}f_{j+1}}-1)\|_{A_p(\mathbb
T)},
$$
and, since $f_{j+1}=0$ on $[0, 2\pi]\setminus I_{j+1}$, we see
that
$$
\|e^{i\lambda S_{j+1}}-e^{i\lambda S_j}\|_{A_p(\mathbb T)}\leq c_{p, 8}
\|(e^{i\lambda \varepsilon_{j+1}f_{j+1}}-1)\|_{A_p(\mathbb T)}.
$$
Therefore (see (58)),
$$
\|e^{i\lambda S_{j+1}}-e^{i\lambda S_j}\|_{A_p(\mathbb T)}\leq c_{p, 9}
2^{-j/q} \Theta_p (\lambda), \qquad \lambda\geq 2, \quad j=0, 1, 2,
\ldots.
$$

   In addition, since
$S_0=\varepsilon_0 f_0=\varepsilon_0 \varphi_0$, we have (see
(56))
$$
\|e^{i\lambda S_0}\|_{A_p(\mathbb T)}=\|e^{i\lambda
\varepsilon_0\varphi_0}\|_{A_p(\mathbb T)}\leq c_{p, 4}\Theta_p(\lambda),
\qquad \lambda\geq 2.
$$

  Thus,
$$
\|e^{i\lambda \varphi}\|_{A_p(\mathbb T)}\leq \|e^{i\lambda
S_0}\|_{A_p(\mathbb T)}+\sum_{j=0}^\infty \|e^{i\lambda
S_{j+1}}-e^{i\lambda S_j}\|_{A_p(\mathbb T)}\leq c_p\Theta_p(\lambda).
$$
The theorem is proved.

\quad

\emph{Remark} 4. The derivative $\varphi'$ of the function
$\varphi$ constructed in Theorem 2 is of bounded variation on
$\mathbb T$. This follows from (37), (38) since for any $m$ the
interval $[0, 2\pi]$ can be partitioned into three intervals so
that the derivative $f_m'$ is monotone on each of them.

\quad

\begin{center}
\textbf{\S~3. Superposition operator}
\end{center}

   In this section we consider the maps of the circle
$\mathbb T$ into itself, i.e., the functions $\varphi : \mathbb
R\rightarrow\mathbb R$ satisfying
$$
\varphi(t+2\pi)=\varphi(t) ~(\mathrm{mod ~2\pi}). \eqno(59)
$$
If such a function is continuous or of class $C^1(\mathbb R)$, or
of class $C^{1, \omega}(\mathbb R)$, then $\varphi$ is
respectively continuous or $C^1$ -smooth, or $C^{1, \omega}$
-smooth map of the circle into itself. If a function $\varphi$
satisfying (59) is continuous and is one-to-one $\textrm{mod}
~2\pi$, then $\varphi$ is a homeomorphism of the circle. If a
homeomorphism of the circle $\varphi$ and its inverse
$\varphi^{-1}$ are $C^1$ -smooth, then $\varphi$ is a
diffeomorphism of the circle. If a diffeomorphism $\varphi$ is
$C^{1,\omega}$ -smooth, then $\varphi$ is called $C^{1, \omega}$
-diffeomorphism (it is easy to verify that in this case the
inverse map $\varphi^{-1}$ is $C^{1, \omega}$ -smooth as well).

  Theorems 3, 4 below follow from Theorems 1, 2
respectively and give their versions for maps of the circle into
itself and integer frequencies.

\quad

\textsc{Theorem 3.} \emph{Let $1\leq p<2$. Let $\varphi$ be a
nonlinear $C^{1, \omega}$ -smooth map of the circle into itself.
Then}
$$
\|e^{in\varphi}\|_{A_p(\mathbb T)}\geq c_p |n|^{1/p}\chi^{-1}(1/|n|),
\qquad n\in\mathbb Z, \quad n\neq 0.
$$

\quad

\textsc{Proof.} We have $\varphi(t+2\pi)=\varphi(t)+2\pi k$. It is
clear that $k\in\mathbb Z$ is independent of $t$. Put
$$
\varphi_0(t)=\varphi(t)-kt.
$$
We obtain $\varphi_0(t+2\pi)=\varphi_0(t)$ and thus $\varphi_0$ is
a real function of class $C^{1, \omega}(\mathbb T)$. The function
$\varphi_0$ is non-constant. It remains to note that
$\|e^{in\varphi}\|_{A_p(\mathbb
T)}=\|e^{in\varphi_0}\|_{A_p(\mathbb T)}$ and apply Theorem 1 to
$\varphi_0$. The theorem is proved.

\quad

  We put $\Theta_1=\Theta$, where $\Theta$ is the function
defined by (13). For $p>1$ the functions $\Theta_p$ are defined by
(26).

\quad

\textsc{Theorem 4.} \emph{Suppose that
$\omega(2\delta)<2\omega(\delta)$ for all sufficiently small
$\delta>0$. There exists a nowhere linear $C^{1, \omega}$
-diffeomorphism $h$ of the circle $\mathbb T$ such that
$$
\|e^{inh}\|_{A_p(\mathbb T)}=O(\Theta_p(|n|)), \qquad
|n|\rightarrow\infty, \quad n\in\mathbb Z,
$$
for all $p, ~1\leq p<2$.}

\quad

\textsc{Proof.} Take the function $\varphi$ from Theorem 2 and put
$$
h(t)=t+\varepsilon\varphi(t),
$$
where $0<\varepsilon\leq 1$ is sufficiently small. We obtain a
nowhere linear diffeomorphism $h$ of the circle $\mathbb T$ of
class $C^{1, \omega}$. It remains to note that
$$
\|e^{inh}\|_{A_p(\mathbb
T)}=\|e^{in\varepsilon\varphi}\|_{A_p(\mathbb T)}\leq
c_p\Theta_p(|n|\varepsilon)\leq c_p\Theta_p(|n|).
$$
The theorem is proved.

\quad

  The corresponding versions of Corollaries 1--3
(as well as of estimate (1)) are obvious. It is also clear that we
have $\|e^{in\varphi}\|_{A_p(\mathbb T)}\simeq
|n|^{\frac{1}{p}-\frac{1}{2}}$ for any nonlinear  $C^{1, 1}$
-smooth map $\varphi : \mathbb T \rightarrow \mathbb T$ (see
Remark 1).

   We shall give now certain natural and obvious
applications of Theorems 3, 4  to the problem on changes of
variable in the spaces $A_p$.

  Let $1<p<2$ and let $\varphi : \mathbb T \rightarrow \mathbb T$
be a continuous map such that for any $f\in A(\mathbb T)$ we have
$f\circ\varphi\in A_p(\mathbb T)$. In this case we say that
$\varphi$ acts from $A$ to $A_p$. Standard arguments (the closed
graph theorem) show that this is the case if and only if the
superposition operator $f\rightarrow f\circ\varphi$ is a bounded
operator from $A(\mathbb T)$ to $A_p(\mathbb T)$. This, in its
turn, is equivalent to the condition
$\|e^{in\varphi}\|_{A_p(\mathbb T)}=O(1), ~n\in\mathbb Z$.

  Using Theorem 3, we see that
if $\omega(\delta)=o(\delta^{p-1})$, then every $C^{1, \omega}$
-smooth map $\varphi$ of the circle into itself that acts from $A$
to $A_p$ is linear. At the same time from Theorem 4 it follows
that if $p>1$, then for any $\varepsilon>0$ there exists a nowhere
linear $C^{1, p-1-\varepsilon}$ -diffeomorphism $h$ of the circle
$\mathbb T$ such that $\|e^{inh}\|_{A_p(\mathbb T)}=O(1)$ and
therefore $h$ acts from $A$ to $A_p$. It is unknown to the author
if one can take here $\varepsilon=0$ (see Remark 3).

   Similarly it is easy to show that the existence of non-trivial
$C^{1, \omega}$ -smooth changes of variable that transfer
functions from $A(\mathbb T)$ to $\bigcap_{p>1} A_p(\mathbb T)$ is
equivalent to the condition that $\omega(\delta)$ decreases to
zero slower than any power, i.e. to the condition
$$
\lim_{\delta\rightarrow+0}\frac{\omega(\delta)}{\delta^\varepsilon}
=\infty \qquad\textrm{for all} ~\varepsilon>0. \eqno(60)
$$
The necessity of this condition follows from Theorem 3. The
sufficiency follows from Theorem 4 (condition (60) implies that
$\Theta_p(y)=O(1), ~y\rightarrow\infty,$ for all $p, ~1<p<2$).
Moreover, if condition (60) holds, then there is a nowhere linear
$C^{1, \omega}$ -diffeomorphism of the circle $h$ such that the
corresponding superposition operator is bounded from $A(\mathbb
T)$ to $A_p(\mathbb T)$ for all $p>1$.

   Concerning this diffeomorphism we note that
the superposition operator $f\rightarrow f\circ h$, which it
generates, is bounded from $A_p(\mathbb T)$ to
$A_{p+\varepsilon}(\mathbb T)$ for any $p, ~1\leq p<2,$ and any
$\varepsilon>0$. Indeed, this superposition operator is a bounded
operator from $A_2(\mathbb T)=L^2(\mathbb T)$ to itself (every
homeomorphism whose inverse satisfies the Lipschitz condition of
order $1$ generates a bounded superposition operator from
$L^2(\mathbb T)$ to itself) and it remains to interpolate $l^p$
between $l^1$ and $l^2$. At the same time we note, that, as we
showed earlier jointly with Olevski\v{\i} [7], every $C^1$ -smooth
map of the circle into itself that for some $p\neq 2$ generates a
bounded superposition operator from $A_p$ to itself is linear.

  In conclusion we recall one open problem that we posed in [9]:
is there a non-trivial map $\varphi : \mathbb T\rightarrow\mathbb
T$ such that
$$
\|e^{in_k\varphi}\|_{A(\mathbb T)}=O(1)
$$
where $n_k, ~k=1, 2, \ldots,$ is some (unbounded) sequence of
integers. As the author showed [9], such a map can not be
absolutely continuous. \footnote{There is a minor inconsistency in
the proof (see [9, Remark]). Everywhere in the proof one should
replace integer $n$ by real $\lambda$.}

\quad

   The results of this paper were partially presented
at the II-d International Symposium on Fourier Series and
Applications, Durso, 2002 (see [16]) and at harmonic analysis
conference HARP 2006 (see [17]).

\quad

   I am grateful to E. A. Gorin for some useful remarks
which improved the paper.

\begin{center} \textbf{References}
\end{center}
\flushleft
\begin{enumerate}

\item A. Beurling, H. Helson, ``Fourier--Stieltjes transforms
    with bounded powers'', \emph{Math. Scand.,} \textbf{1}
    (1953), 120-126.

\item J.-P. Kahane, \emph{S\'erie de Fourier absolument
    convergentes}, Springer-Verlag, Berlin--Heidelberg--New
    York, 1970.

\item J.-P. Kahane, ``Quatre le\c cons sur les
    hom\'eomorphismes du circle et les s\'eries de Fourier'',
    in: \emph{Topics in Modern Harmonic Analysis,} Vol. II,
    Ist. Naz. Alta Mat. Francesco Severi, Roma, 1983, 955-990.

\item  Z. L. Leibenson, ``On the ring of functions with
    absolutely convergent Fourier series'', \emph{Uspehi
    Matem. Nauk}, \textbf{9}:3(61) (1954), 157-162 (in
    Russian).

\item  J.-P. Kahane, ``Sur certaines classes de s\'eries de
    Fourier absolument convergentes'', \emph{J. de
    Math\'ematiques Pures et Appliqu\'ees}, \textbf{35}:3
    (1956), 249-259.

\item  L. Alpar, ``Sur une classe partiquli\`ere de s\'eries
    de Fourier \`a certaines puissances absolument
    convergentes'', \emph{Studia Sci. Math. Hungarica},
    \textbf{3} (1968), 279-286.

\item  V. Lebedev, A. Olevski\v{\i}, ``$C^1$ changes of
    variable: Beurling -- Helson type theorem and H\"ormander
    conjecture on Fourier multipliers'', \emph{Geometric and
    Functional Analysis (GAFA)}, \textbf{4}:2 (1994), 213-235.

\item  A. Zygmund, \emph{Trigonometric series}, vol. 1,
    Cambridge Univ. Press, New York, 1959.

\item  V. V. Lebedev, ``Diffeomorphisms of the circle and the
    Beurling--Helson theorem'', \emph{Functional analysis and
    its applications}, \textbf{36}:1(2002), 25-29.

\item  M. N. Leblanc, ``Sur la r\'eciproque de l'in\'egalit\'e
    de Carlson'', \emph{C.R. Acad. Sc. Paris, S\'erie A,}
    \textbf{267} (1968), 332-334.

\item  N. K. Bari, \emph{A treatise on trigonometric series},
    vols. I, II, Pergamon Press, Oxford, 1964.

\item  E. M. Stein, G. Weiss, \emph{Introduction to Fourier
    analysis on Euclidean spaces}, Princeton Math. Ser., vol.
    32, Princeton Univ. Press, Princeton, NJ, 1971.

\item  M. Plancherel, G. P\'olya, ``Fonctions enti\`eres et
    int\'egrales de Fourier multiples. II'', \emph{Comment.
    Math. Helv.}, \textbf{10}:2 (1937), 110--163.

\item E. M. Stein, \emph{Singular integrals and
    differentiability properties of functions}, Princeton
    Math. Ser., vol. 30, Princeton Univ. Press, Princeton, NJ,
    1970.

\item  J. L. Rubio de Francia, ``A Littlewood--Paley
    inequality for arbitrary intervals'', \emph{Rev. Mat.
    Iberoam.,} \textbf{1}:2 (1985), 1-14.

\item  V. V. Lebedev, ``Quantitative estimates in the
    Beurling--Helson theorem'', \emph{II Int. Symp. Fourier
    Series and Applications, Durso, May 27 -- June 2, 2002.
    Abstracts}, Rostov--na--Donu, 2002, 33-34.

\item  V. Lebedev, ``Rate of growth in Beurling -- Helson
    theorem'', \emph{Harmonic Analysis and Related Problems,
    HARP 2006, Zaros, Crete, Greece, June 19 - 23, 2006},
    Abstracts.

\end{enumerate}

\quad

\qquad\textsc{V. V. Lebedev}\\
\qquad Dept. of Mathematical Analysis\\
\qquad Moscow State Institute of Electronics\\
\qquad and Mathematics (Technical University)\\
\qquad e-mail: \emph {lebedevhome@gmail.com}

\end{document}